%% file: OrdinalConeOverleaf.tex
\documentclass[preprint]{elsarticle}
\usepackage[utf8]{inputenc}
\usepackage{pdfpages}
\usepackage[hidelinks]{hyperref}
\usepackage{url}
\usepackage{amsthm}
\usepackage{thmtools}
\usepackage{amssymb}
\usepackage{amsmath}
\usepackage{float}
\usepackage{pdflscape}
\usepackage{natbib}
\usepackage{todonotes}
\usepackage{changes}
\definechangesauthor[name=MS, color=green!50!black]{MS}
\definechangesauthor[name=KK, color=pink]{KK}
\definechangesauthor[name=JS, color=blue]{JS}

\usepackage{dsfont}

\usepackage[labelformat=simple]{subcaption}
\usepackage{tikz}
\usetikzlibrary{patterns}
\usetikzlibrary{matrix}
\usetikzlibrary{graphs,quotes}\usepackage{pgfplots}
\pgfplotsset{compat=1.15,
	every axis/.append style={
		axis lines=center,
		xlabel style={anchor=south west},
		ylabel style={anchor=south west},
		zlabel style={anchor=south west},
		tick align=outside,}
}
\usepgfplotslibrary{patchplots}
\pgfdeclareverticalshading{darkGray}{100bp}
{color(0bp)=(black!45); color(32bp)=(black!45); color(38bp)=(black!40);color(45bp)=(black!30); color(100bp)=(black!25)}

\pgfdeclareverticalshading{darkGray2}{100bp}
{color(0bp)=(black!45); color(40bp)=(black!45); color(45bp)=(black!40);color(50bp)=(black!30); color(100bp)=(black!25)}

\pgfdeclareverticalshading{midGray}{100bp}
{color(0bp)=(black!30); color(30bp)=(black!35);color(35bp)=(black!30);color(40bp)=(black!25); color(47bp)=(black!15); color(100bp)=(black!10)}

\pgfdeclareverticalshading{lightGray}{100bp}
{color(0bp)=(black!20); color(28bp)=(black!15);color(35bp)=(black!15); color(40bp)=(black!10); color(45bp)=(black!5); color(100bp)=(black!5)}

\pgfdeclareverticalshading{lightGray2}{100bp}
{color(0bp)=(black!20); color(35bp)=(black!15);color(40bp)=(black!15); color(45bp)=(black!10); color(50bp)=(black!5); color(100bp)=(white)}

\usepackage{cancel}

\usepackage[ruled,vlined,linesnumbered]{algorithm2e}
\SetKwProg{Fn}{Function}{}{}
\DontPrintSemicolon
\SetKwFor{For}{for}{}{}
\SetKwFor{While}{while}{}{}

\SetCommentSty{mycommfont}

\usepackage{rotating} %ermöglicht es figures zu drehen (\begin{sidewaysfigure})
\usepackage{booktabs} % \hline wird in Tabelle schöner (mehr Abstand oben und unten plaziert, man muss aber \midrule oder \toprule oder \bottomrule verwenden) oder durch \addlinespace[vspace] mehr abstand zwischen Zeilen in Tabellen erzeugen
\usepackage{mathtools}

\newtheorem{defi}{Definition}
\newtheorem{theorem}[defi]{Theorem}
\newtheorem{lem}[defi]{Lemma}

\newtheorem{rem}[defi]{Remark}

\newtheorem{example}[defi]{Example}

\newcommand {\Z}{\mathds{Z}}
\newcommand {\Zgeq}{\mathds{Z}_\geq}

\newcommand {\R}{\mathds{R}}

\DeclareMathOperator{\argmin}{argmin}

\DeclareMathOperator{\sort}{sort}
\DeclareMathOperator{\ordmin}{ordinally\text{ }minimize}
\DeclareMathOperator{\interior}{int} %interior
\DeclareMathOperator{\rank}{rank}

\newcommand{\zulM}{X} %zulässige Menge allgemein
\DeclareMathOperator{\eff}{eff}
\newcommand{\effM}{X_{\eff}} % effiziente Menge allgemein
\DeclareMathOperator{\nd}{nd}
\newcommand{\dS}{S} % diskrete Menge
\newcommand{\zulTM}{X} %zulässige (Teil-)Menge allgemein (der Potenzmenge der diskreten Menge)
\newcommand{\xone}{x'} %erster zulässiger Punkt
\newcommand{\xtwo}{\hat{x}} % zweiter zulässiger Punkt
\newcommand{\xthree}{\bar{x}} % dritter zulässiger Punkt
\newcommand{\outNdD}{Y_{\nd}} % Menge der nicht-dominierten Outcome Vectoren

\makeatletter
\newcommand*\bigcdot{\mathpalette\bigcdot@{1}}
\newcommand*\bigcdot@[2]{\mathbin{\vcenter{\hbox{\scalebox{#2}{$\m@th#1\bullet$}}}}}
\makeatother
\newcommand{\pareto}{P} %Zeichen für den Pareto Kegel ohne 0
\newcommand{\vR}{\nu} % numerical representation
\newcommand{\VR}{\mathcal{V}} % set of all numerical representation
\newcommand{\preceqqnu}{\preceqq} % weakly ordinally-dominates Zeichen
\newcommand{\preceqnu}{\preceq} % ordinally-dominates Zeichen
\newcommand{\cone}{C} %Zeichen für einen cone
\newcommand{\tailC}{\cone_{\leqslant_{\lsum}}} %Zeichen für den Tail Cone
\newcommand{\tailA}{A_{\leqslant_{\lsum}}} % Matrix für Tail-Dominance über Halbräume
\newcommand{\tailB}{B_{\leqslant_{\lsum}}} % Matrix für Tail-Dominance über Extremstrahlen
\newcommand{\headA}{A_{\geqslant_{\fsum}}} % Matrix für Head-Dominance über Halbräume
\newcommand{\headB}{B_{\geqslant_{\fsum}}} % Matrix für Head-Dominance über Extremstrahlen
\newcommand{\headC}{C_{\geqslant_{\fsum}}} %Zeichen für den Head Cone
\DeclareMathOperator{\hcone}{hcone} %Symbol für die Funktion die aus einer Matrix den cone der Halfspaces berechnet
\DeclareMathOperator{\vcone}{vcone} %Symbol für die Funktion die aus einer Matrix den cone der vectors (extreme rays) berechnet
\DeclareMathOperator{\OCOP}{OCOP}
\DeclareMathOperator{\OOP}{OOP}
\DeclareMathOperator{\TOP}{TOP}

\newcommand{\lsum}{t} %symbol für summe über last elements
\newcommand{\fsum}{h} %symbol für summe über first elements

\newcommand{\preceqq}{\;\raisebox{-0.12em}{\scalebox{1.15}{\ensuremath{\mathrel{\substack{\prec\\[-.15em]=}}}}}\;}

\makeatletter
\addtotheorempostfoothook{%
   \global\everypar{{\setbox\z@\lastbox}\global\everypar{}}%
}
\makeatother

\usepackage{numcompress}\biboptions{authoryear}
\defcitealias{geovelo}{\slshape geovelo}

\begin{document}

\title{Ordinal Optimization Through Multi-objective Reformulation}
\author[1]{Kathrin Klamroth}
\ead{klamroth@math.uni-wuppertal.de}

\author[1]{Michael Stiglmayr}
\ead{stiglmayr@math.uni-wuppertal.de}

\author[1]{Julia Sudhoff\corref{cor1}}
\ead{sudhoff@math.uni-wuppertal.de}

\cortext[cor1]{Corresponding author}
\address[1]{University of Wuppertal, School of Mathematics and Natural Sciences, IMACM, Gaußstr.~20, 42119 Wuppertal, Germany}

\begin{abstract}
    We analyze combinatorial optimization problems with ordinal, i.e., non-additive, objective functions that assign categories (like good, medium and bad) rather than cost coefficients to the elements of feasible solutions. We review  different optimality concepts  
	%In combinatorial optimization already exist different definitions of optimality 
	for ordinal optimization problems and discuss their similarities and differences. We then focus on two prevalent optimality concepts that are shown to be equivalent. % In this article, we investigate three different definitions of ordinal optimality. We show that all three definitions are equivalent under certain conditions and that two of them are equivalent in general. 
	Our main result is a bijective linear transformation that transforms ordinal optimization problems to associated standard multi-objective optimization problems  with binary cost coefficients. Since this transformation preserves all properties of the underlying problem, problem-specific solution methods remain applicable. A prominent example is dynamic programming and Bellman's principle of optimality, that can be applied, e.g., to ordinal shortest path and ordinal knapsack problems. We extend our results to multi-objective optimization problems that combine ordinal and real-valued objective functions.
\end{abstract}

\maketitle

% \tableofcontents

%%%%%%%%%%%%%%%%%%%%%%%%%%%%%%%%%%%%%%%%
\section{Introduction}
\label{sec:intro}
%%%%%%%%%%%%%%%%%%%%%%%%%%%%%%%%%%%%%%%%
Ordinal objective functions occur whenever it is only possible to represent the quality of an element by a ordered category and not by a numerical value. As an example, consider the problem of finding optimized routes for cyclists in a road network: While edges may be associated with different categories like asphalt, gravel or sand---or, when related to safety considerations, very safe (there is a bicycle path), neutral (a quiet road) or unsafe (a main road without bicycle path)---such categories do not immediately translate into monetary or cost values.
Bi-objective shortest path problems with route safety criteria are addressed, for example, in the web application \citetalias{geovelo} and in the associated publications \cite{kerg:anef:2021,sauv:sear:2010}. In these references, only two categories are considered (safe or unsafe edges), and the safety criterion is translated into a cost function that evaluates the total length of unsafe route segments.
In contrast, an ordinal shortest path problem is investigated in \cite{SCHAFER20}. A major difficulty when considering ordinal objective functions is that ``optimality'' may be defined in many different ways. \cite{SCHAFER20} suggests an optimality concept that is based on sorted category vectors. A similar concept is used in \cite{klamroth2021multiobjective}, where matroid optimization problems with one real-valued and one ordinal objective function are investigated. A different perspective is proposed in \cite{SCHAFER:knapsack} who define ordinal optimality for knapsack problems on the basis of numerical representations for the categories.
 
 In this paper, we consider general combinatorial optimization problems and provide a new cone-based interpretation of the optimiality concept for ordinal objectives suggested in \cite{SCHAFER:knapsack}. In particular, we interrelate ordinal optimality with the classical concept of Pareto optimality for an associated multi-objective optimization problem. For a general introduction to multi-objective optimization we refer to \cite{Ehrg05}. Since the underlying transformation of the objective function is linear and bijective and hence preserves the  combinatorial structure of the respective problems, our results immediately lead to efficient solution strategies as, for example, dynamic programming for shortest path problems. 
The respective transformations are based on a representation of dominance relations by cones.

The paper is organized as follows. In Section~\ref{sec:basics} we review three optimality concepts for ordinal optimization based on  \cite{SCHAFER:knapsack}. These results and optimality concepts %of \cite{SCHAFER:knapsack}
are then extended and re-interpreted as a special case of cone-optimality \citep[see, e.g.,][]{engau2007domination} in Section~\ref{sec:cone}. %Therefore, we define the corresponding cones and the multi-objective optimization problem. 
A detailed analysis of the properties of the corresponding ordering cones then leads to a linear transformation of ordinal optimization problems to associated multi-objective optimization problems with the Pareto cone defining dominance and optimality. This transformation is used in Section~\ref{sec:solution} to formulate a general algorithm for solving ordinal optimization problems. Furthermore, we investigate the relation between the definition of ordinal optimality and the weight space decomposition of the associated multi-objective problem. In Section~\ref{sec:extensions} we extend our results to more general problem types with additional real-valued objective functions. We conclude in Section~\ref{sec:conclusion} with a summary and an outlook on future research.

%%%%%%%%%%%%%%%%%%%%%%%%%%%%%%%%%%%%%%%%
%alter Titel \section{Ordinal Optimality Concepts}
\section{Single-objective Ordinal Optimization}
\label{sec:basics}

\subsection{Problem Definition}\label{subsec:problemDef}
%------------------------------------------------------------------------------------
%\textcolor{yellow!60!black}{ordinal objective function + OOP(ordinal optimization problem) + categories}\\
We consider combinatorial optimization problems with an ordinal objective function. In general, an \emph{ordinal optimization problem} \eqref{eq:OOP} can be formulated as
\begin{equation}\label{eq:OOP}\tag{OOP}
	\begin{array}{rl}
		``\ordmin" & o(x)\\
		\text{s.\,t.} & x\in \zulTM,
	\end{array}
\end{equation}
where $\zulTM$ is the set of feasible solutions. We assume that $\zulTM$ is a subset of the power set of a discrete set $\dS$, i.e., $\zulTM\subseteq2^\dS$.  Every element of $\dS$ is assigned to one of $K$ ordered categories. This assignment is encoded by a mapping $o:\dS\to \mathcal{C}$ with $\mathcal{C}=\{\eta_1,\ldots,\eta_K\}$. We assume that category \(\eta_i\) with \(i\in\{1,\ldots,K-1\}\) is strictly preferred over category \(\eta_{i+1}\), written as $\eta_i \prec \eta_{i+1}$. The objective function of a feasible solution $x=\{e_1,\dots,e_n\}$ is given by the \emph{ordinal vector} $o(x)=\sort(o(e_1),\dots,o(e_n))$, where the operator $\sort(\,)$ means that the components of $o(x)$ are sorted w.r.t.\  non-decreasing preferences, i.e., $o_1(x)\preceqq o_2(x) \preceqq \cdots \preceqq o_n(x)$. Note that different feasible solutions may have different numbers of elements, and hence the length of the ordinal vector $o(x)$ may vary for different $x\in\zulTM$.

%------------------------------------------------------------------------------------------------------------------------------
%\textcolor{yellow!60!black}{Def: counting function}\\
Instead of using the un-aggregated, ordered ordinal vector \(o(x)\) one can count the number of elements in a feasible solution per category.  Accordingly, we use the \emph{counting vector} $c:\zulTM \to\Z_\geqq^K$ with $\Z_\geqq^K\coloneqq \{y\in\Z^K:y_i\geq0 \text{ for all }i=1,\dots,K\}$. Thereby, the \(i\)-th component of $c(x)$ equals the number of elements in $x$ which are in category $\eta_i$, i.\,e., 
$c_i(x)=\vert\{e\in x \colon o(e)=\eta_i\}\vert$. 
Obviously, there is a one to one correspondence between the vectors $o\in\R^n$ and $c\in\R^K$, since the ordinal vector $o$ can be determined from a given counting vector $c$ by $o_i(x)=\eta_j$ with $j=\argmin\{j\in\{1,\dots,K\}:i\leq\sum_{l=1}^j c_l(x)\}$.  Note again that the number of elements $n$ of a feasible solution, and hence the length of the ordinal vectors $o\in\R^n$, may vary while the number of categories $K$ and therefore the length of the counting vectors $c\in\R^K$ is fixed. Hence, we get the following formulation of an \emph{ordinal counting optimization problem} \eqref{eq:OCOP}
\begin{equation}\label{eq:OCOP}\tag{OCOP}
	\begin{array}{rl}
		``\ordmin" & c(x)\\
		\text{s.\,t.} & x\in \zulTM,
	\end{array}
\end{equation}
which will be shown to be equivalent to problem~\eqref{eq:OOP} for an appropriate definition of ``ordinal minimization".

In the following, we also consider an \emph{incremental tail counting vector} $\tilde{c}\in\R^K$ that counts, in its $i$-th component, the number of elements of a feasible solution $x$ which are  in category $\eta_i$ or worse, i.e., $\tilde{c}_i(x)=\vert \{e\in x\colon \eta_i\preceqq o(e)\} \vert=\sum_{j=i}^K c_j(x)$.
In particular, the total number of elements of a solution $x$ is given in the first component of $\tilde{c}$, i.e., \(|x|= \tilde{c}_1 (x)=\sum_{i=1}^{K} c_i(x)\).
%------------------------------------------------------------------------------------------------------------------------------

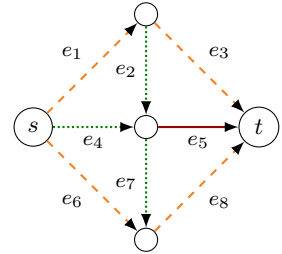
\begin{figure}[htb]
	\centering
		\small
		\hspace*{\fill}
		\parbox{0.58\textwidth}{
		\begin{tabular}{c|ccc}
            \toprule
			&$o$&$c$&$\tilde{c}$\\
			\midrule 
			\rule{0pt}{25pt}$x^1=\{e_1,e_2,e_5\}$&$\begin{pmatrix} \eta_1\\ \eta_2\\ \eta_3 \end{pmatrix}$&$\begin{pmatrix} 1\\ 1\\ 1 \end{pmatrix}$&$\begin{pmatrix} 3\\ 2\\ 1 \end{pmatrix}$\\
			\rule{0pt}{25pt}$x^2=\{e_4,e_5\}$&$\begin{pmatrix} \eta_1\\ \eta_3 \end{pmatrix}$&$\begin{pmatrix} 1\\ 0\\ 1 \end{pmatrix}$&$\begin{pmatrix} 2\\ 1\\ 1 \end{pmatrix}$\\ 
			\rule{0pt}{25pt}$x^3=\{e_1,e_3\}$&$\begin{pmatrix} \eta_2\\ \eta_2 \end{pmatrix}$&$\begin{pmatrix} 0\\ 2\\ 0 \end{pmatrix}$&$\begin{pmatrix} 2\\ 2\\ 0 \end{pmatrix}$\\
			\rule{0pt}{25pt}$x^4=\{e_6,e_8\}$&$\begin{pmatrix} \eta_2\\ \eta_2 \end{pmatrix}$&$\begin{pmatrix} 0\\ 2\\ 0 \end{pmatrix}$&$\begin{pmatrix} 2\\ 2\\ 0 \end{pmatrix}$\\
			\rule{0pt}{25pt}$x^5=\{e_4,e_7,e_8\}$&$\begin{pmatrix} \eta_1\\ \eta_1\\ \eta_2 \end{pmatrix}$&$\begin{pmatrix} 2\\ 1\\ 0 \end{pmatrix}$&$\begin{pmatrix} 3\\ 1\\ 0 \end{pmatrix}$\\ 
			\rule{0pt}{30pt}$x^6=\{e_1,e_2,e_7,e_8\}$&$\begin{pmatrix} \eta_1\\ \eta_1\\ \eta_2\\ \eta_2 \end{pmatrix}$&$\begin{pmatrix} 2\\ 2\\ 0 \end{pmatrix}$&$\begin{pmatrix} 4\\ 2\\ 0 \end{pmatrix}$\\ \bottomrule
		\end{tabular}}
		\hspace*{\fill}
		\hspace*{\fill}
        \parbox{0.3\textwidth}{
		\input{Bilder/bspDirGraph}
		}
	\caption{Instance of an ordinal shortest path problem. A dotted-green edge is in the best category $\eta_1$, a dashed-orange edge is in category $\eta_2$ and a solid-red edge is in the worst category $\eta_3$. All possible $s-t$ paths $x^i$, $i=1,\dots,6$ and their respective objective function vectors are given. }
	\label{fig:shortestPath1}
\end{figure}

%\textcolor{yellow!60!black}{Beispiel}\\
As an example, consider the shortest path problem shown in Figure~\ref{fig:shortestPath1} together with the outcome vectors $o(x)$ (for problem \eqref{eq:OOP}) and $c(x)$ (for problem \eqref{eq:OCOP}) for all feasible solutions $x\in\zulTM$. In addition, the incremental tail counting vector $\tilde{c}(x)$ is given for all $x\in\zulTM$.

%----------------------------------------------------------------------------------------------------------------------------

%%%%%%%%%%%%%%%%%%%%%%%%%%%%%%%%%%%%%%%%%%%%%%%%%%%%%%%%%%%%%%%%%%%%%%%%%%%%%%%%%%%%%%%%%%%%%%%%%%%%%%%%%%%%%%%%%%%%%%%%%%%%%%%%%%

\subsection{Optimality Concepts for Ordinal Objective Functions}\label{subsec:optConcepts}
%\textcolor{yellow!60!black}{Question: what does minimization mean? Idea: for every value I choose for each category, a solution is always non-dominated}\\
In the following, we review three different concepts of optimality for ordinal optimization. All of them try to answer the question what minimization could mean for the problems \eqref{eq:OOP} and \eqref{eq:OCOP}. The first concept is to use \emph{numerical representations} that assign a numerical value to every category such that the order of the categories is respected. In this context, a numerical representation respects the order of the categories whenever the numerical value of a better category is strictly smaller than the numerical value of a less preferred category.  If we take the sum over all numerical values of a vector $o(\xone)$ for a feasible solution $\xone$, we obtain a unique numerical value that can be compared to the corresponding numerical value of another feasible solution $\xtwo$. %We call a feasible solution $\xone$ efficient if there is no numerical representation such that there is an other feasible solution $\xtwo$ with smaller value.
A feasible solution $\xone$ is called \emph{efficient} if there is no other feasible solution $\xtwo$ which is better w.r.t.\ \emph{all} numerical representations.

The second concept is to maximize the number of elements in the good categories, and the third concept is to minimize the number of elements in the bad categories. 
After the formal introduction of these three optimality concepts, we investigate their interrelation.

%------------------------------------------------------------------------------------------------------------------------------
\paragraph{Optimality by Numerical Representations}
%\textcolor{yellow!60!black}{Def: Numerical representation}\\
We consider the problems \eqref{eq:OOP} and \eqref{eq:OCOP}. The concept of \emph{optimality by numerical representation} for ordinal objectives as introduced in \cite{SCHAFER:knapsack} is based on a previous and more general work of \cite{FISHBURN1999359}. It assigns an order preserving numerical value to each category. %A formal definition of a \emph{numerical representation} which is similar to the definition in \cite{SCHAFER:knapsack} is provided below. 
%For a survey on numerical representations of \kk{general} preference structures see \cite{FISHBURN1999359}. 
Following \cite{SCHAFER:knapsack}, we call a function $\vR:\mathcal{C}\to\Zgeq$ a \emph{numerical representation} if 
\begin{equation*}
	\eta_i \prec \eta_j \iff \vR(\eta_i)<\vR(\eta_j)\text{ for all }i,j\in\{1,\dots,K\}.
\end{equation*}
 Note that we assume strictly ordered categories, i.e., there are no categories that are indifferent. As a consequence, we do not allow $\vR(\eta_i)=\vR(\eta_j)$ for $i\neq j$ since this would make two different categories indistinguishable in the numerical representation. Let $\VR$ denote the set of all numerical representations for a given number of categories $K$.

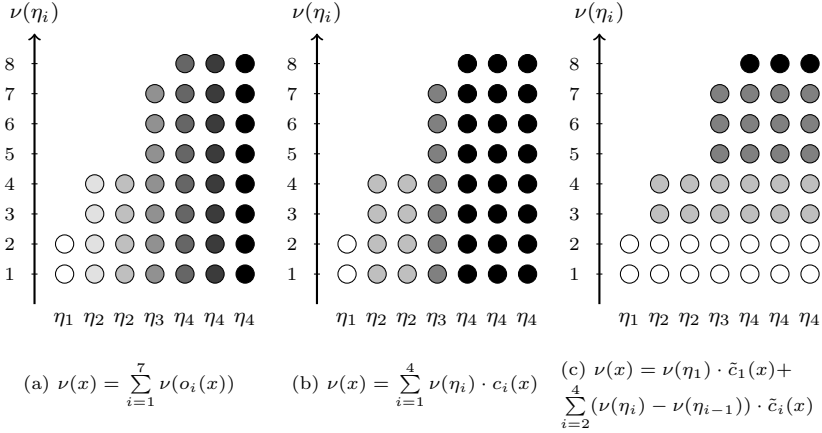
\begin{figure}
	\centering
	\subcaptionbox{$\vR(x)=\sum\limits_{i=1}^7 \vR(o_i(x))$ \label{fig:numRep1}}
	[.3\linewidth]{\input{Bilder/numerRepraesentation1-neu}}
	\subcaptionbox{$\vR(x)=\sum\limits_{i=1}^4 \vR(\eta_i)\cdot c_i(x)$\label{fig:numRep2}}
	[.3\linewidth]{\input{Bilder/numerRepraesentation2}}
	\subcaptionbox{$\vR(x)=\vR(\eta_1)\cdot\tilde{c}_1(x)+$\\$\sum\limits_{i=2}^4(\vR(\eta_i)-\vR(\eta_{i-1}))\cdot \tilde{c}_i(x)$\label{fig:numRep3}}
	[.3\linewidth]{\input{Bilder/numerRepraesentation3}}
	\caption{Consider an example with $n=7$ and $K=4$. Different ways to compute the numerical value of a feasible solution $x$ with $o(x)=(\eta_1,\eta_2,\eta_2,\eta_3,\eta_4,\eta_4,\eta_4)^\top$, $c(x)=(1,2,1,3)^\top$, $\vR(\eta_1)=2$, $\vR(\eta_2)=4$, $\vR(\eta_3)=7$ and $\vR(\eta_4)=8$ are illustrated. The different colors represent the summands and visualize the different slicing strategies. }
	\label{fig:numRepr}
\end{figure}

%------------------------------------------------------------------------------------------------------------------------------
%\textcolor{yellow!60!black}{Def: Numerical Value of a feasible solution}\\
For a given numerical representation $\vR$, we define the numerical value of a feasible solution $x=\{e_1,\dots,e_n\}\in\zulTM$ w.r.t. $\vR$ \citep[cf.][]{SCHAFER:knapsack} as 
\begin{align*}
	\vR(x)&\coloneqq\sum_{i=1}^n\vR(o(e_i)) =\sum_{i=1}^K\vR(\eta_i)\cdot c_i(x).
\end{align*}
The numerical value $\vR(x)$ of a feasible solution $x\in\zulTM$ can be evaluated in different ways by re-arranging the terms and using the counting vector $c$ or the incremental tail counting vector $\tilde{c}$, respectively:
\begin{align*}
	\vR(x)&=\sum_{i=1}^K\vR(\eta_i)\cdot c_i(x)\\
	&=\sum_{i=1}^{K-1} \vR(\eta_i)\biggl(\sum_{j=i}^K c_j(x)-\sum_{j=i+1}^K c_j(x)\biggr)+\vR(\eta_K)\,c_K(x)\\
	&=\vR(\eta_1)\cdot\sum_{i=1}^K c_i(x)+\sum_{i=2}^K \bigl(\vR(\eta_i)-\vR(\eta_{i-1})\bigr)\cdot \sum_{j=i}^K c_j(x)\\
	&=\vR(\eta_1)\cdot\tilde{c}_1(x)+\sum_{i=2}^K\bigl(\vR(\eta_i)-\vR(\eta_{i-1})\bigr)\cdot \tilde{c}_i(x).
\end{align*}
An illustration for different ways to evaluate $\vR(x)$ is given in Figure~\ref{fig:numRepr}. %\todo{ordered median problem}

%-----------------------------------------------
%\textcolor{yellow!60!black}{Beispiel}\\
\begin{example}\label{ex:numRep}
We apply the concept of numerical representations to the shortest path problem given in Figure~\ref{fig:shortestPath1}. %to motivate the resulting optimality concept. 
As a motivation, suppose that there are two decision makers $A$ and $B$ who have to select a most preferred path. They would agree, for example, that path $x^1$ is worse than path $x^2$, because $x^1$ has a dashed-orange edge more than $x^2$ and, other than that, their outcome vectors are the same. But they do not agree on the question whether $x^2$ or $x^5$ is preferred, because decision maker $A$ chooses the numerical representation $\vR_A(\eta_1)=1$, $\vR_A(\eta_2)=2$ and $\vR_A(\eta_3)=5$, while decision maker $B$ chooses $\vR_B(\eta_1)=2$, $\vR_B(\eta_2)=3$ and $\vR_B(\eta_3)=4$. Therefore, decision maker $A$ would prefer path $x^5$ because $\vR_A(x^5)=4<6=\vR_A(x^2)$ while decision maker $B$ would prefer path $x^2$ because $\vR_B(x^2)=6<7=\vR_B(x^5)$. Hence, the path $x^2$ does \emph{not} ordinally dominate the path $x^5$, i.e., $x^2$ is \emph{not} better than $x^5$ for \emph{all} numerical representations. Similarly, the path $x^5$ does not ordinally dominate the path $x^2$ since also $x^5$ is not better than $x^2$ for all numerical representations. %Hence, we call a path $x^*$ ordinally efficient if there is no other path $x$ which is better, or as good as, $x^*$ for all numerical representations. This means in our example that the paths $x^2$ and $x^5$ both are ordinally efficient and that path $x^1$ is ordinally dominated by $x^2$. Furthermore, the paths $x^3$ and $x^4$ are also ordinally efficient.
\end{example}
%------------------------------------------------------------------------------------------------------------------------------
%\textcolor{yellow!60!black}{Efficiency Dominance}
\begin{defi}[cf.~\citealp{SCHAFER:knapsack}]
Let $\xone,\xtwo\in\zulTM$ be feasible solutions. Then,
\begin{enumerate}
	\item $\xone$ \emph{weakly ordinally dominates} $\xtwo$, $o(\xone)$ \emph{weakly ordinally dominates} $o(\xtwo)$ and $c(\xone)$ \emph{weakly ordinally dominates} $c(\xtwo)$, denoted by $\xone\preceqqnu \xtwo$, $o(\xone)\preceqqnu o(\xtwo)$, $c(\xone)\preceqqnu c(\xtwo)$, respectively, if and only if for \emph{every} $\vR\in\VR$, it holds that $\vR(\xone)\leq \vR(\xtwo)$.
	\item $\xone$ \emph{ordinally dominates} $\xtwo$, $o(\xone)$ \emph{ordinally dominates} $o(\xtwo)$ and $c(\xone)$ \emph{ordi\-nally dominates} $c(\xtwo)$, denoted by $\xone\preceqnu \xtwo$, $o(\xone)\preceqnu o(\xtwo)$, $c(\xone)\preceqnu c(\xtwo)$, respectively, if and only if  $\xone$ weakly ordinally dominates $\xtwo$ and there exists $\vR^*\in\VR$ such that $\vR^*(\xone)< \vR^*(\xtwo)$.
	\item $x^*\in\zulM$ is called \emph{ordinally efficient}, if there does not exist an $x\in\zulTM$ such that $x\preceqnu x^*$.
	\item $o(x^*)$ and $c(x^*)$ are called \emph{ordinally non-dominated outcome vectors} of Problem~\eqref{eq:OOP} and \eqref{eq:OCOP}, respectively, if $x^*$ is ordinally efficient.
\end{enumerate}
\end{defi}

%------------------------------------------------------------------------------------------------------------------------------
\paragraph{Optimality by Maximization of Elements in Good Categories}
%\textcolor{yellow!60!black}{Optimalität mit Summe über erste Elemente}\\
Another optimality concept in ordinal optimization is to maximize the number of elements in good categories. The intuition behind this concept is that solutions with many good elements are to be preferred over solutions with few good elements. The drawback, however, is that this concept rewards solutions with larger numbers of elements as long as these are in (relatively) good categories, which may not be wanted in practice. This optimality concept is defined only for the problem \eqref{eq:OCOP}, as we need the counting vector $c$ for its definition.
\begin{defi}\label{def:headsum}
We say $\xone$ \emph{weakly head-dominates} $\xtwo$, denoted by $\xone\geqq_{\fsum} \xtwo$ or $c(\xone)\geqq_{\fsum} c(\xtwo$), if and only if
\begin{equation}\label{eq:sumFirstElem}% \tag{$headSum$}
	\sum_{i=1}^j c_i(\xone) \geq \sum_{i=1}^j c_i(\xtwo) \;\text{ for all } j=1,\ldots ,K.
\end{equation}
Furthermore, $\xone$ \emph{head-dominates} $\xtwo$, denoted by $\xone\geqslant_{\fsum} \xtwo$ or $c(\xone)\geqslant_{\fsum} c(\xtwo$), if and only if \eqref{eq:sumFirstElem} holds and $c(\xone)\neq c(\xtwo)$. Moreover, $x^*\in\zulTM$ is called \emph{head-efficient} if there is no $x\in\zulTM$ such that $x\geqslant_{\fsum}x^*$. The corresponding outcome vector $c(x^*)$ is called \emph{head-non-dominated}. %\todo{JS: wollen wir wirklich auch  $o(x^*)$ hier als head-non-dominated bezeichnen? Brauchen wir eigentlich nicht, oder? (Analog für tail...)}
\end{defi}

%------------------------------------------------------------------------------------------------------------------------------
\paragraph{Optimality by Minimization of Elements in Bad Categories}
%\textcolor{yellow!60!black}{Optimalität mit Summe über letzte Elemente}\\
The drawback that longer solutions may be preferred over shorter solutions, as long as the elements are in good categories, can be avoided by taking the converse perspective, i.e., when minimizing the number of elements in the bad categories. Again, this optimality concept is defined only for the problem \eqref{eq:OCOP}.
\begin{defi}\label{def:tailsum}
	We say $\xone$ \emph{weakly tail-dominates} $\xtwo$, denoted by $\xone\leqq_{\lsum} \xtwo$ or $c(\xone)\leqq_{\lsum} c(\xtwo$), if and only if
\begin{equation}\label{eq:sumLastElem}%\tag{$tailSum$}
	\tilde{c}_j(\xone)=\sum_{i=j}^K c_i(\xone) \leq \sum_{i=j}^K c_i(\xtwo)=\tilde{c}_j(\xtwo) \;\text{ for all }j=1,\ldots ,K.
\end{equation}
Again, $\xone$ \emph{tail-dominates} $\xtwo$, denoted by $\xone\leqslant_{\lsum} \xtwo$ or $c(\xone)\leqslant_{\lsum} c(\xtwo$), if and only if \eqref{eq:sumLastElem} holds and $c(\xone)\neq c(\xtwo)$. Moreover, $x^*\in\zulTM$ is called \emph{tail-efficient} if there is no $x\in\zulTM$ such that $x\leqslant_{\lsum}x^*$.  The corresponding outcome vector $c(x^*)$ is called \emph{tail-non-dominated}.
\end{defi}

\begin{rem}
Note that head-dominance as well as tail-dominance are equivalently defined on the feasible set $\zulTM\subseteq 2^\dS$ and on its image set $c(X)\subseteq\R^K$. The definitions immediately extends to the complete $\R^K$.
\end{rem}

%------------------------------------------------------------------------------------------------------------------------------
\subsection{Properties of and Interrelations between  Optimality Concepts for Ordinal Optimization}\label{sub:relationOC}
%\textcolor{yellow!60!black}{Remark other Definitions of Optimality are possible}\\

In addition to the concepts described above, there are further ways to define efficiency. This has been done, for example, in \cite{SCHAFER20} for the ordinal shortest path problem. Their definition has the disadvantage that \emph{Bellman's principle of optimality} \citep[see][]{bellman57dynamic} does not hold in general, i.e.,\ not every subpath of an efficient path is necessarily efficient w.r.t.\ this optimality concept. The definition of head-optimality has the same disadvantage, see Remark~\ref{rem:headsum} below for more details. In contrast, the definitions of ordinal optimality and tail-optimality can be proven to be equivalent. Moreover, they are compliant with Bellman's principle of optimality. Note that, for the special case of a knapsack problem, this was shown in \cite{SCHAFER:knapsack}.

%------------------------------------------------------------------------------------------------------------------------------
%\textcolor{yellow!60!black}{Rem: Definition is consistent, i.e., gibt es nur Elemente einer Kategorie dominiert die Lösung mit den wenigsten Elementen (weglassen?)}\\
%\begin{rem}
%	The definition of ordinally efficient solutions is consistent in the sense that, if there are two feasible solutions of different cardinality but with all elements from the same category, then the solution with less elements ordinally dominates the solution with more elements. 
%\end{rem} 
%This observation leads to the following result.
%\todo[inline]{KK: Vielleicht sogar beides: Erst die Erläuterung, dann die formale Aussage.}
\begin{lem}\label{lem:Teilmenge}
	For feasible solutions $\bar{x}=\{\bar{e}_1,\dots,\bar{e}_n\},x'=\{e'_1,\dots,e'_m\}\in\zulTM$ with 
	$c_i(\bar{x})\leq c_i(x')$ for all $i=1,\dots,K$
	%$\{o(\bar{x}_1),\dots,o(\bar{x}_n)\}\subseteq \{o(x'_1),\dots,o(x'_m)\}$ 
	and $n<m$ 
	it holds $\bar{x}\preceqnu x'$.
\end{lem}
%\todo[inline]{KK Wir sollten uns eine konsistente Notation für die Kardinalität einer Lösung überlegen. k ist eher ungünstig (assoziiert mit Kategorien). Vielleicht s,t,...? Oder l,n,m,..?.}
\begin{proof}
    First let $\vR\in\VR$ be an arbitrary numerical representation. %Wlog we assume that $o(\bar{x}_i)=o(x'_i)$ for all $i=1,\dots,n$. 
	Then $\vR(\bar{x})=%\sum_{i=1}^n \vR(o(\bar{x}_i))=\sum_{i=1}^n \vR(o(x'_i))\leq \sum_{i=1}^m \vR(o(x'_i))
	\sum_{i=1}^K \vR(\eta_i)\cdot c_i(\bar{x}) \leq \sum_{i=1}^K \vR(\eta_i)\cdot c_i(x')=\vR(x')$, i.e.,\ $\bar{x}\preceqqnu x'$.
	Note that since $n<m$, and due to the above assumptions, there must exist a category $\eta_j$, $j\in\{1,\dots,K\}$ such that $c_j(\bar{x})<c_j(x')$. Therefore, there exists a numerical representation $\vR^*$ such that $\vR^*(\eta_j)>0$ and thus $\vR^*(\bar{x})<\vR^*(x')$, which concludes the proof.
\end{proof} %\todo{JS: den ersten Satz direkt for den letzten schieben? Die Aussage wird vorher nicht gebraucht.}
Note that the condition of Lemma \ref{lem:Teilmenge} is always satisfied if $\bar{x},x'\in\zulTM$ with $\bar{x}\subsetneq x'$. Thus,  in Example~\ref{ex:numRep} the path $x^2$ is always preferred over the path $x^1$. This is also the case for tail-efficiency, but not for head-efficiency, see Remark~\ref{rem:headsum} below.
In many application contexts this is meaningful property, since adding additional elements to a solution (no matter from which category) does generally not improve the solution quality. As an example, we refer again to paths representing bicycle routes as in the app \citetalias{geovelo}, where we are not interested in routes that are unnecessarily long.

%------------------------------------------------------------------------------------------------------------------------------
%\textcolor{yellow!60!black}{Lem: dominance Relation is a preorder}\\
%Lem 3.7 in \cite{SCHAFER:knapsack}\\
\begin{lem}[\citealt{SCHAFER:knapsack}]\label{lem:relation}
	The ordinal dominance relation $\preceqqnu$ defined on the feasible set $\zulTM$ is a preorder, i.e., it is reflexive and transitive.
\end{lem}
\begin{proof}
	Let $\xone,\xtwo,\xthree\in\zulTM$. Obviously, $\vR(\xone)\leq\vR(\xone)$ holds for every $\vR\in\VR$. Hence, the relation $\preceqqnu$ is reflexive. If $\vR(\xone)\leq\vR(\xtwo)$ and $\vR(\xtwo)\leq\vR(\xthree)$ for every $\vR\in\VR$ it follows that $\vR(\xone)\leq\vR(\xthree)$ for every $\vR\in\VR$ by definition and therefore, we have shown transitivity.
\end{proof}
Note that the ordinal dominance relation $\preceqqnu$ is in general not antisymmetric on the feasible set $\zulTM$ since two different feasible solutions may have the same number of elements in each category like the paths $x^3$ and $x^4$ in Example~\ref{ex:numRep}.
%\todo[inline]{KK unklar formuliert. Besser: Beispiel angeben?}

The following two results show that (weak) ordinal dominance and (weak) tail-dominance are actually equivalent on the feasible set $X$.

%------------------------------------------------------------------------------------------------------------------------------
%\textcolor{yellow!60!black}{Lem: Bezug dominance relation und counting function (weak efficiency)}\\
%Lem 4.1 in \cite{SCHAFER:knapsack}\\
\begin{lem}[\citealt{SCHAFER:knapsack}]\label{lem:weakEff}
	Let $\xone,\xtwo\in\zulTM$ be two feasible solutions. Then $\xone$ weakly ordinally dominates $\xtwo$, i.e., $\xone\preceqqnu \xtwo$ if and only if $\xone\leqq_{\lsum} \xtwo$.
\end{lem}
%\todo[inline]{KK Was haltet Ihr davon, mit dieser Definition für eine ordinale Ordnung zu starten? Das könnte man auch gut motiviern (es ist immer besser, weniger teure Elemente zu haben). Dann könnte man auf die Äquivalenz zu der Definition über $\nu$ verweisen und den Beweis weglassen. Den Beweis zu wiederholen ist denke ich nicht gut; direkt fällt mir keine deutliche Vereinfachung ein. Das sollten wir nochmal diskutieren, also nicht direkt etwas ändern.}
\begin{proof}
	The proof is a simplified variant of the proof in \cite{SCHAFER:knapsack}. Note, that they consider maximization problems while we consider minimization problems.
	
	First we show by contradiction that $\xone\preceqqnu \xtwo$ implies $\xone\leqq_{\lsum} \xtwo$. Let $\xone,\xtwo\in\zulTM$ and let $j^*\in\{1,\dots,K\}$ with 
	\begin{equation*}
		\sum_{i=j^*}^{K} c_i(\xone) > \sum_{i=j^*}^{K} c_i(\xtwo).
	\end{equation*}
	The idea of the proof is to make the bad categories $\eta_{j^*},\dots,\eta_K$ very expensive, such that an element of this category can not be replaced by elements of the lower categories. Hence, we define the numerical representation 
	\begin{equation*}%\label{eq:numRep}
		\vR(\eta_i)=\begin{cases}
			i,&\text{ if }i<j^*\\
			i+2\, \vert \xtwo\vert\, K,&\text{ if }i\geq j^*.
		\end{cases}
	\end{equation*}

	This implies
	\begin{align*}
		\vR(\xone)&\geq 2\,|\xtwo|\,K \cdot \sum_{i=j^*}^K c_i(\xone) \geq 2\,|\xtwo|\,K \cdot \biggl(1+\sum_{i=j^*}^K c_i(\xtwo) \biggr)\\
		& > |\xtwo|\,K + 2\,|\xtwo|\,K\cdot\sum_{i=j^*}^K c_i(\xtwo) \geq \sum_{i=1}^K i\, c_i(\xtwo) + 2\,|\xtwo|\,K\cdot\sum_{i=j^*}^K c_i(\xtwo)\\
		& = \sum_{i=1}^{j^*-1} i\, c_i(\xtwo) +\sum_{i=j^*}^K i\, c_i(\xtwo)+\sum_{i=j^*}^K 2\,|\xtwo|K\cdot c_i(\xtwo)= \vR(\xtwo) .
	\end{align*} 

	For the other direction we use the reformulation of $\vR(x)$, which is visualized in Figure~\ref{fig:numRep3}. It follows that for any $\vR\in\VR$
	\begin{align*}
		\vR(\xone)&=\vR(\eta_1)\,\tilde{c}_1(\xone)+\sum_{i=2}^K \bigl(\vR(\eta_i)-\vR(\eta_{i-1})\bigr)\,\tilde{c}_i(\xone)\\
		&\leq \vR(\eta_1)\,\tilde{c}_1(\xtwo)+\sum_{i=2}^K \bigl(\vR(\eta_i)-\vR(\eta_{i-1})\bigr)\,\tilde{c}_i(\xtwo) = \vR(\xtwo).
	\end{align*}
	The inequality holds because of the assumption $\tilde{c}_j(\xone)\leq \tilde{c}_j(\xtwo)$ for all $j=1,\dots,K$ and $\vR(\eta_i)-\vR(\eta_{i-1})>0$ for all $\vR\in\VR$ and $i=2,\dots,K$.
	Hence, we have shown that $\sum_{i=j}^K c_i(\xone) \leq \sum_{i=j}^K c_i(\xtwo)$ for all $j=1,\ldots ,K$ implies $\xone\preceqqnu \xtwo$, which concludes the proof.
\end{proof}

%------------------------------------------------------------------------------------------------------------------------------
%\textcolor{yellow!60!black}{Lem: Bezug dominance relation und counting function (efficiency)}\\
\begin{lem}\label{lem:eff}
	Let $\xone,\xtwo\in\zulTM$ be two feasible solutions. Then $\xone$ ordinally dominates $\xtwo$, i.e., $\xone\preceqnu \xtwo$ if and only if $\xone\leqslant_{\lsum} \xtwo$.
\end{lem}
\begin{proof}
    We first show that $\xone\preceqnu\xtwo$ implies $\xone\leqslant_{\lsum} \xtwo$.
	If $\xone$ ordinally dominates $\xtwo$, then $\xone\preceqqnu \xtwo$ which implies $\xone\leqq_{\lsum}\xtwo$ due to Lemma \ref{lem:weakEff}. It remains to show that $c(\xone)\neq c(\xtwo)$. As $\xone$ ordinally dominates $\xtwo$ it holds that there is a numerical representation $\vR^*$ such that $\vR^*(\xone)<\vR^*(\xtwo)$. Hence, 
	\begin{align*}
		&0<\vR^*(\xtwo)-\vR^*(\xone)\\
		\iff\,&0<\sum_{i=1}^K \vR^*(\eta_i) \bigl(c_i(\xtwo)-c_i(\xone)\bigr).
	\end{align*}
	 Since $\vR(\eta_i)>0$ for all $i=1,\dots,K$, it holds $c(\xone)\neq c(\xtwo)$. Consequently, we have shown that when $\xone$ ordinally dominates $\xtwo$, then $\xone\leqq_{\lsum}\xtwo$ holds for all $j=1,\ldots ,K$ and $c(\xone)\neq c(\xtwo)$.
	 
	 For the other direction it is sufficient to show that $c(\xone)\neq c(\xtwo)$ implies that there exists a numerical representation $\vR^*\in\VR$ such that $\vR^*(\xone)<\vR^*(\xtwo)$. Let $j^*$ be the largest category such that $c_{j^*}(\xone)\neq c_{j^*}(\xtwo)$. This implies 
	 \begin{equation*}
	 	\sum_{i=j^*}^{K} c_i(\xone) < \sum_{i=j^*}^{K} c_i(\xtwo).
	 \end{equation*} Now the result follows analogously to the proof of Lemma~\ref{lem:weakEff} with exchanged roles of $\xone$ and $\xtwo$.
\end{proof}

%------------------------------------------------------------------------------------------------------------------------------
%\textcolor{yellow!60!black}{Rem: Lemma is not true, if we consider the sum over the first elements}\\
\begin{rem}\label{rem:headsum}
Lemma \ref{lem:weakEff} (and thus also Lemma \ref{lem:eff}) does not hold in general for the relation $\geqslant_{\fsum}$. %In other words, it is not true in general that for two feasible solutions $\xone,\xtwo\in\zulTM$,   $\xone$ weakly ordinally dominates $\xtwo$, i.e., $\xone\preceqqnu \xtwo$ if and only if \eqref{eq:sumFirstElem} holds.
		As a counter example, consider the paths $x^1$ and $x^2$ with $c(x^1)=(1,1,1)^\top$ and $c(x^2)=(1,0,1)^\top$ from Figure~\ref{fig:shortestPath1}. Obviously, $x^1$ head-dominates $x^2$. But for every numerical representation it follows $\vR(x^1)>\vR(x^2)$, which contradicts $x^1\preceqqnu x^2$.
\end{rem}

Note that the crucial point in the counter example given in Remark~\ref{rem:headsum} is the different cardinality of the solutions, \(|x^1|\neq|x^2|\). For ordinal optimization problems with fixed cardinality, for which matroid optimization problems as studied in \cite{klamroth2021multiobjective} are an example, it can be shown that head- and tail-dominance are equivalent.

%------------------------------------------------------------------------------------------------------------------------------
%\textcolor{yellow!60!black}{Lem: Both summerizing strategies lead to the same solution if all feasible solutions have the same length (weglassen?)}\\
\begin{lem}
	If all feasible solutions have the same cardinality, i.e., if \(|\xone|=|\xtwo|\) for all \(\xone,\xtwo\in \zulTM\),
% 	\begin{align}\label{eq:sameLength}
% 		&\sum_{i=1}^K c_i(\xone)=\sum_{i=1}^K c_i(\xtwo) \quad \text{for all }\xone,\xtwo\in \zulTM
% 	\end{align} 
	then head- and tail-dominance as defined in Definitions~\ref{def:headsum} and \ref{def:tailsum}, respectively, are equivalent.
\end{lem}

\begin{proof}
	First assume that $\xone$ head-dominates $\xtwo$, i.e., inequality~\eqref{eq:sumFirstElem} is satisfied. We show that then $\xone$ also tail-dominates $\xtwo$, i.e., inequality~\eqref{eq:sumLastElem} holds. Towards this end, let $j\in\{2,\dots,K\}$. Then
	\begin{align*}
		\sum_{i=j}^K c_i(\xone)&=\sum_{i=1}^K c_i(\xone)-\sum_{i=1}^{j-1} c_i(\xone)\overset{\eqref{eq:sumFirstElem}}{\leq}\sum_{i=1}^K c_i(\xone)-\sum_{i=1}^{j-1} c_i(\xtwo)&\\
		&\overset{|\xone|=|\xtwo|}{=}\sum_{i=1}^K c_i(\xtwo)-\sum_{i=1}^{j-1} c_i(\xtwo)=\sum_{i=j}^K c_i(\xtwo),&
	\end{align*}
	which implies \eqref{eq:sumLastElem}.
	
	Now let $\xone$ tail-dominate $\xtwo$, i.e., \eqref{eq:sumLastElem} is satisfied. We show that then also $\xone$ head-dominates $\xtwo$, i.e., \eqref{eq:sumFirstElem} holds. Hence, let $j\in\{1,\dots,K-1\}$. Then
	\begin{align*}
		\sum_{i=1}^j c_i(\xone)&=\sum_{i=1}^K c_i(\xone)-\sum_{i=j+1}^K c_i(\xone)\overset{\eqref{eq:sumLastElem}}{\geq}\sum_{i=1}^K c_i(\xone)-\sum_{i=j+1}^K c_i(\xtwo)&\\
		&\overset{|\xone|=|\xtwo|}{=}\sum_{i=1}^K c_i(\xtwo)-\sum_{i=j+1}^K c_i(\xtwo)=\sum_{i=1}^j c_i(\xtwo),&
	\end{align*}
	which implies \eqref{eq:sumFirstElem}.
\end{proof}

%------------------------------------------------------------------------------------------------------------------------------
%\textcolor{yellow!60!black}{Remark: sum relation is (strict) partial order and compatible with addition/ scalar multiplication}\\
\begin{lem}
	The relation $\leqq_{\lsum}$ is a partial order on $\R^K$, i.e., it is reflexive, transitive and antisymmetric. Moreover, the relation $\leqslant_{\lsum}$ is a strict partial order on $\R^K$, i.e., it is irreflexive and transitive.%\deleted[id=MS]{, on $\R_\geqq^K\coloneqq\{y\in\R^K:y_i\geq0\text{ for all }i=1,\dots,K\}$.} 
\end{lem}
%\todo[inline]{JS: Beweis so lassen? bei den Lemmas beziehen sich die relationen auf X und nicht $\R^K$}
\begin{proof}
    %The relation $\leqq_{\lsum}$ is equivalent to the relation $\preceqqnu$ due to Lemma~\ref{lem:weakEff},  and hence it is reflexive and transitive due to Lemma~\ref{lem:relation}.
	Let $u\in\R^K$. Then $u\leqq_{\lsum} u$, i.e.,\ $\leqq_{\lsum}$ is reflexive. Furthermore, for $u,v,w\in\R^K$ such that $u\leqq_{\lsum} v$ and $v\leqq_{\lsum} w$ it follows that $\sum_{i=j}^K u_i \leq \sum_{i=j}^K v_i \leq \sum_{i=j}^K w_i$ for all $j=1,\dots,K$, i.e.,\ $u\leqq_{\lsum} w$ which means $\leqq_{\lsum}$ is transitive. 
	To show that the relation $\leqq_{\lsum}$ is antisymmetric, consider two vectors $u,v\in\R^K$ with $u\leqq_{\lsum} v$ and $v\leqq_{\lsum} u$. Then
	 $\sum_{i=j}^K u_i = \sum_{i=j}^K v_i$ for all $j=1,\dots,K$. This implies that $u=v$  and hence $\leqq_{\lsum}$ is antisymmetric. Therefore, $\leqq_{\lsum}$ is a partial order.
	
	Now consider the relation $\leqslant_{\lsum}$.
	Since $u\nleqslant_{\lsum} u$ for all $u\in\R^K$, it holds that $\leqslant_{\lsum}$ is irreflexive.	It remains to show that $\leqslant_{\lsum}$ is transitive. Towards this end, consider three vectors $u,v,w\in\R^K$ such that $u\leqslant_{\lsum} v$ and $v\leqslant_{\lsum} w$. This implies $\sum_{i=j}^K u_i \leq \sum_{i=j}^K v_i \leq \sum_{i=j}^K w_i$ for all $j=1,\dots,K$, and there exist indices $s,t\in\{1,\dots,K\}$ such that $\sum_{i=s}^K u_i < \sum_{i=s}^K v_i$ and $\sum_{i=t}^K v_i < \sum_{i=t}^K w_i$. Hence, we can conclude that $\sum_{i=j}^K u_i \leq \sum_{i=j}^K w_i$ for all $j=1,\dots,K$ and $u\neq w$, i.e.,\ $u\leqslant_{\lsum} w$. Consequently, we have shown that $\leqslant_{\lsum}$ is irreflexive and transitive which concludes the proof.
\end{proof}

As a consequence of the discussion in Sections~\ref{subsec:optConcepts} and \ref{sub:relationOC}, 
we focus in the following on the ordinal optimization problem \eqref{eq:OOP} w.r.t.\ optimality by numerical representation, or equivalently, on the ordinal counting optimization problem \eqref{eq:OCOP} w.r.t.\ the concept of tail-dominance. %Since the latter turns out to be easy to evaluate and implement, we often refer to problem \eqref{eq:OCOP} and tail-dominance when talking about \emph{ordinal optimality} in the following.

%\todo[inline]{JS: Genau genommen steht bei den Problemen (OOP) und (OCOP) nur "ordinally minimize", dh wie minimiert wird ist zu dem Zeitpunkt nicht spezifiziert. (OOP) kann man nur mit der ordinalen dominance ($\min_\prec$) betrachten, aber für (OCOP) kann man doch an dieser Stelle schon sagen, dass man "Ordinally minimize" equivalent ersetzen kann durch $\min_\prec$ oder gleichwertig durch $\min_{\leqslant_{\lsum}}$. Damit können wir im Grunde auch an dieser Stelle auch bei (OCOP) von einem Ordinal Problem und ordinal optimality sprechen, weil wir das Problem mit dieser Optimalität definieren können - aus meiner Sicht ist damit das Problem mit den Begrifflichkeiten insofern gelöst, dass wir (OCOP) mit ordinal bezeichnen können, gleichzeitig müssen wir jetzt immer genau angeben, welches optimalitäts-Konzept wir gerade meinen (Problem angeben reicht nicht aus), denn im Folgenden setzen wir erstmal die tail dominance in Bezug zu den Kegeln. Ich bin mittlerweile auch der Meinung, dass der Abschnitt hier bleiben kann und nicht in die Einleitung von Kapitel 3 gehört.}
%\todo[inline]{KK: Ich habe oben den Nachsatz wieder eingefügt - so jetzt ok? Ich würde ungern eine künstliche Unetrscheidung einführen bei Problemen, die äquivalent sind.}

%%%%%%%%%%%%%%%%%%%%%%%%%%%%%%%%%%%%%%%%%%%%%%%%%%%%%%%%%%%%%%%%%%%%%%%%%%%%%%%%%%%%%%%%%%%%%%%%%%%%%%%%%%%%%%%%%%%%%%%%%%%%%%%%%%%

%alter Titel: \section{Ordering Cones and Optimality}
\section{Ordinal Optimality versus Pareto Optimality: An Interpretation Based on Ordering Cones}
\label{sec:cone}
%------------------------------------------------------------------------------------------------------------------------------

%\textcolor{yellow!60!black}{Example: größer (gleich)Ordnungen auf dem reelen Vektorraum}
A prominent example of an order relation in the context of multi-objective optimization is the  \emph{compo\-nent-wise order} or \emph{Pareto order}, see, e.g., \cite{Ehrg05}. For two vectors $u,v\in\R^K$, we write
\begin{align}\label{eq:componentwise}
%	u\leqq v &\,:\Longleftrightarrow\, u_i \leq v_i,\; \quad i=1,\ldots ,K,\\
	u\leqslant v &\,:\Longleftrightarrow\, u_i\leq v_i,\; \quad i=1,\ldots,K \;\;\text{and}\;\; u\neq v,%\\
%	u<v &\,:\Longleftrightarrow\, u_i<v_i, \;\quad i=1,\ldots,K,
\end{align}
and say that $u$ \emph{Pareto dominates} $v$ if and only if $u\leqslant v$. The associated weak and strict component-wise orders, respectively, are defined by
\begin{align}
	u\leqq v &\,:\Longleftrightarrow\, u_i \leq v_i,\; \quad i=1,\ldots ,K, \label{eq:weakcomponentwise}\\
%	u\leqslant v &\,:\Longleftrightarrow\, u_i\leq v_i,\; \quad i=1,\ldots,K \;\;\text{and}\;\; u\neq v,%\\
	u<v &\,:\Longleftrightarrow\, u_i<v_i, \;\quad i=1,\ldots,K.\label{eq:strictcomponentwise}
\end{align}
Note that $\leqq$ defines a partial order in $\R^K$ while $\leqslant$ as well as $<$ define strict partial orders in $\R^K$. 

In the following, basic properties of orders and their relation to cones are discussed and illustrated at the example of the Pareto order. This leads to a novel perspective on ordinal optimality in comparison to Pareto optimality.

%alter Titel: \subsection{Introduction of Cones}
\subsection{Orders and Cones}
\label{subsec:intro}
Orders and cones are closely related. 
The following review of basic concepts  relevant in our context is  on \cite{Ehrg05,engau2007domination,ziegler95}.

A \emph{cone} in $\R^K$ is a subset $\cone \subseteq\R^K$ such that $\lambda \, u  \in \cone$ for all $u\in\cone$ and for all $\lambda\in\R$, $\lambda>0$.
A cone $\cone\in\R^K$ is called \emph{pointed} if $u\in \cone$ implies that $(-u)\not\in \cone$ for all $u\neq 0$. 
%------------------------------------------------------------------------------------------------------------------------------
%\textcolor{yellow!60!black}{Definition Polyhedral cone}\\
%(Quelle Dissertation, oder auch \cite{ziegler95})\\
Moreover, a cone $\cone\subseteq\R^K$ is called a \emph{polyhedral cone} if there exists a matrix $A\in\R^{m\times K}\setminus\{0\}$ such that $\cone=\hcone(A)\coloneqq\{y\in\R^K:A\,y\geqq0\}$. The rows of the matrix $A$ are normal vectors of hyperplanes, and thus a polyhedral cone can be seen as the finite intersection of $m$ (closed and linear) halfspaces. Polyhedral cones can also be described by their extreme rays. This property is an immediate consequence of the well-known Weyl-Minkowski Theorem:
%\textcolor{yellow!60!black}{Weyl-Minkowski Theorem}\\
\begin{theorem}[Weyl-Minkowski-Theorem, cf.\ \citealt{ziegler95}] \label{th:W-M}
	A cone $\cone\subseteq \R^K$ is finitely generated by $n$ vectors in $\R^K$, i.e.,
	\begin{equation*}
		\cone= \vcone(B)\coloneqq\bigl\{B\lambda:\lambda\in\R^n,\lambda\geqq 0\bigr\}\;\text{ for some }B\in\R^{K\times n}
	\end{equation*}
	if and only if it is a finite intersection of $m$ halfspaces in $\R^{K}$, i.e.,
	\begin{equation*}
		\cone=\hcone(A)=\bigl\{y\in\R^K:A\,y\geqq0\bigr\}\; \text{ for some }A\in\R^{m\times K}.
	\end{equation*}
\end{theorem}
%This result shows that a polyhedral cone can equivalently be described by its extreme rays. 
%\kk{Note that in the two descriptions given in Theorem~\ref{th:W-M}, every row of $A$ corresponds to the normal vector of a hyperplane defining $\cone$, while} every column in $B$ corresponds to a point on an extreme ray of $\cone$. \todo[inline]{JS: Das steht unmittelbar vor dem Theorem schon fast genauso da - streichen? Übergang zu den dualen Kegeln?}

%------------------------------------------------------------------------------------------------------------------------------
%\textcolor{yellow!60!black}{Def: Dualer Kegel}\\
Now let $\cone\subset \R^K$ be a cone. % set, in particular not necessarily a cone. 
Then the sets 
\begin{align*}
	&\cone^*\coloneqq\{d\in\R^K:d^\top c\geq0\text{ for all }c\in\cone\} \\
	&\cone^*_s\coloneqq\{d\in\R^K:d^\top c>0\text{ for all }c\in\cone\setminus\{0\}\}
\end{align*}
are called the \emph{dual cone} and the \emph{strict dual cone} of $\cone$, respectively.

%----------------------------------------------------------------------------------------------------

Every cone $\cone\in\R^K$ induces a binary (ordering) relation $\mathcal{R}\subseteq\R^K\times\R^K$ 
by defining that $(u,v)\in\mathcal{R}$ if and only if $(v-u)\in \cone$. 
Depending on the context, we also write $u\mathcal{R}v$ whenever $(u,v)\in\mathcal{R}$ (as, for example, in the case of the binary relations defined in \eqref{eq:componentwise}, \eqref{eq:weakcomponentwise} and \eqref{eq:strictcomponentwise} above).
Binary relations that are induced by cones are always  compatible with scalar multiplication, i.e., $(u,v)\in\mathcal{R}$ implies $(\lambda\,  u,\lambda\, v)\in\mathcal{R}$ for all $u,v\in\R^K$ and $\lambda>0$. Moreover, they are compatible with addition, i.e., $(u,v)\in\mathcal{R}$ implies $(u+w,v+w)\in\mathcal{R}$ for all $u,v,w\in\R^K$.

Conversely, binary (ordering) relations that are compatible with scalar multiplication  induce cones that represent the respective relation. The following result will be particularly useful in our context.

%------------------------------------------------------------------------------------------------------------------------------
%\textcolor{yellow!60!black}{Lem: cone induced by binary relation which is compatible with scalar multiplication and addition (Prop 1.14 + Thm 1.17 Ehrgott)}\\
\begin{lem}[see, e.g., \citealt{Ehrg05}]\label{lem:ordcone}
	Let $\mathcal{R}\subseteq\R^K\times\R^K$ be a binary relation on $\R^K$ which is compatible with scalar multiplication. Then $\cone_\mathcal{R}\coloneqq\{(v-u)\in\R^K:(u,v)\in\mathcal{R}\}$ is a cone, and $\cone_\mathcal{R}$ induces the binary relation $\mathcal{R}$. If $\mathcal{R}$ is additionally compatible with addition, then the following statements hold:
	\begin{enumerate}
		\item $0\in\cone_\mathcal{R}$ if and only if $\mathcal{R}$ is reflexive.
		\item $\cone_\mathcal{R}$ is pointed if and only if $\mathcal{R}$ is antisymmetric.
		\item $\cone_\mathcal{R}$ is convex if and only if $\mathcal{R}$ is transitive.
	\end{enumerate}
\end{lem}

%-----------------------------------------------------------------------------------------------------------------------------

%
It is easy to see that all three %component-wise 
orders \eqref{eq:componentwise}, \eqref{eq:weakcomponentwise} and \eqref{eq:strictcomponentwise} %$\leqslant$, $\leqq$ and $<$, 
are compatible with scalar multiplication and addition. The component-wise order \eqref{eq:componentwise} induces the \emph{Pareto cone} given by $\pareto\coloneqq\R_\geqslant^K=\{(v-u)\in\R^K : u\leqslant v\}=\{y\in\R^K:y\geqslant0\}$. %, which is the cone induced by the strict partial order $\geqslant$.  
Similarly, the weak component-wise order \eqref{eq:weakcomponentwise} induces the %\textcolor{yellow!60!black}{Def: Pareto Cone}\\
\emph{closure of the Pareto cone given by} $\R_\geqq^K=\{y\in\R^K: y\geqq 0 %y_i\geq0\text{ for all }i=1,\dots,K
\}=\pareto\cup \{0\}$, which is a proper, pointed and convex cone (see Lemma~\ref{lem:ordcone}). Moreover, it is a polyhedral cone that is defined by the identity matrix, i.e., $\pareto\cup\{0\}=\hcone(I)=\vcone(I)$ (where $I$ is the $K\times K$ identity matrix). Note also that the Pareto cone is self dual, i.e., $\pareto^*=\pareto$.

Lemma~\ref{lem:ordcone} implies that binary relations $\mathcal{R}$ that are compatible with scalar multiplications can be equivalently represented by associated cones $\cone_{\mathcal{R}}$. This interrelation was used, among others, in \cite{engau2007domination} to define the concept of \emph{cone-efficiency} (or \emph{$\cone_{\mathcal{R}}$-efficiency}). %that is, under appropriate assumptions, equivalent to the concept of $\mathcal{R}$-efficiency as introduced above.
%Towards this end, we refer to the concept of cone optimality as, for example, suggested by \cite{engau2007domination}.
For a general introduction to ordering cones in the context of vector optimization see, e.g., \cite{tammer03theory,jahn11vector}.
\begin{defi}[c.f.~\citealt{engau2007domination}]\label{def:coneOpt}
	Let $Y\subset\R^K$ be a nonempty set and let $\cone_\mathcal{R}\subset\R^K$ be a cone induced by a strict partial order 
	$\mathcal{R}\subset\R^K\times \R^K$ (i.e., $\mathcal{R}$ is irreflexive and transitive). Then the sets
	\begin{align*}
		N(Y,\cone_\mathcal{R})&\coloneqq\{y\in Y: (y-\cone_\mathcal{R})\cap Y=\varnothing\} \\
		N_w(Y,\cone_\mathcal{R})&\coloneqq\{y\in Y: (y-\interior(\cone_\mathcal{R}))\cap Y=\varnothing\}
	\end{align*}
	are called the \emph{$\cone_\mathcal{R}$-non-dominated} and the \emph{weakly $\cone_\mathcal{R}$-non-dominated} set of $Y$, %with respect to $\cone_\mathcal{R}$, 
	respectively. The corresponding pre-images $x\in\zulTM$ are called \emph{$\cone_\mathcal{R}$-efficient} and \emph{weakly $\cone_\mathcal{R}$-efficient}, respectively.
	Thereby, $\interior(\cone_\mathcal{R})$ denotes the interior of $\cone_\mathcal{R}$. 
	
	Furthermore, we say that $u$ \emph{$\cone_\mathcal{R}$-dominates} $v$ if 
	$u\in (v-\cone_\mathcal{R})$, and that $u$ \emph{weakly $\cone_\mathcal{R}$}-dominates $v$ if $u\in (v-\interior(\cone_\mathcal{R}))$.
	%$u\mathcal{R}v$. 
\end{defi} 
%\todo[inline]{JS: $\cone_\mathcal{R}$-efficient doch definieren? Insbesondere nutzen wir später Pareto-efficient}

Note that when $\mathcal{R}$ is the component-wise order (or Pareto order) defined in \eqref{eq:componentwise}, then $\cone_\mathcal{R}=\pareto$ is the Pareto cone in $\R^K$, $N(Y,\pareto)$ is the non-dominated set (or \emph{Pareto set}) of $Y$, and $N_w(Y,\pareto)$ is the weakly non-dominated set of $Y$, see, e.g., \cite{Ehrg05}.

%------------------------------------------------------------------------------------------------------------------------------

\subsection{The Ordinal Cone} 
\label{subsec:ordinalcone}

In the following we show that tail-dominance, represented by the binary relation $\leqslant_{\lsum}$, is induced by a polyhedral cone in $\R^K$. We will refer to this cone as the \emph{ordinal cone}. As a first step towards this goal, we prove that $\leqslant_{\lsum}$ is compatible with scalar multiplication and addition. %i.e., $u\leqslant_{\lsum}v$ implies $\lambda \, u\leqslant_{\lsum}\lambda \, v$ for all $u,v\in\R^K_{\geqq}$ and $\lambda>0$, and $u\leqslant_{\lsum}v$ implies $u+w\leqslant_{\lsum}v+w$ for all $u,v,w\in\R^K_{\geqq}$. 
As a second step, the associated ordinal cone is constructed and analyzed. Afterwards, we can reinterpret problem~\eqref{eq:OCOP} based on cone optimality.

\begin{lem}\label{lem:compatible}
	The relation $\leqslant_{\lsum}$ is compatible with scalar multiplication and with  addition.
\end{lem}
\begin{proof}
To show that $\leqslant_{\lsum}$ is compatible with scalar multiplication, let $\lambda>0$ and $u,v\in\R^K$ with %\todo{KK: Das gilt in ganz $\R^K$ und nicht nur in $\R^K_{\geqq}$}
$u\leqslant_{\lsum} v$. It follows that $\lambda \sum_{i=j}^K u_i\leq \lambda \sum_{i=j}^K v_i$ which implies $ \sum_{i=j}^K \lambda \, u_i\leq \sum_{i=j}^K \lambda \, v_i$ for all $j=1,\dots,K$. Furthermore, it holds that $\lambda \, u\neq \lambda \, v$, and hence $\lambda\, u\leqslant_{\lsum} \lambda\, v$, which implies that $\leqslant_{\lsum}$ is compatible with scalar multiplication.

It remains to show that $\leqslant_{\lsum}$ is also compatible with addition. Let $u,v,w\in\R^K$ with $u\leqslant_{\lsum} v$, i.e.,\ $\sum_{i=j}^K u_i\leq\sum_{i=j}^K v_i$ for all $j=1,\dots,K$ and $u\neq v$. This implies that $\sum_{i=j}^K (u_i+w_i)\leq\sum_{i=j}^K (v_i+w_i)$ for all $j=1,\dots,K$ and $(u+w)\neq (v+w)$, i.e., $u+w\leqslant_{\lsum}v+w$. Hence we have proven the compatibility with addition.
\end{proof}

It can be proven analogously that $\leqq_{\lsum}$ is also compatible with scalar multiplication and addition.

From Lemma~\ref{lem:ordcone} and Lemma~\ref{lem:compatible} we can conclude that the cone $\tailC\coloneqq\{(v-u)\in\R^K:u\leqslant_{\lsum}v\}$ induced by the strict partial order $\leqslant_{\lsum}$ is pointed, convex and it does not contain $0$. We call this cone the \emph{ordinal cone} to emphasize that $\tailC$ equivalently represents ordinal dominance and show that its closure, the cone  $\tailC\cup\{0\}$, is a polyhedral cone that can be described as the intersection of $K$ halfspaces.

%------------------------------------------------------------------------------------------------------------------------------
%\textcolor{yellow!60!black}{Theorem: Der Kegel über die Facetten entspricht dem Kegel induziert durch die binary relation}\\
\begin{theorem}\label{thm:A}
	The closure of the ordinal cone is a polyhedral cone. %\todo{KK: Ich finde es einfacher, die Aussagen für $\tailC\cup\{0\}$ zu zeigen und habe das so geändert.} 
	In particular, it holds that $\tailC\cup\{0\}=\hcone(\tailA)$ 
	with $\tailA\in\R^{K\times K}$ given by %\todo{KK: Die Bezeichnung $A$, $B$ ist nicht eindeutig, daher hier und im Folgenden spezifiziert (vgl.\ z.B.\ Thms.~12 und 22). Bitte nochmal genau durchsuchen und idealerweise einen latex Befehl einführen.}
	\begin{equation*}
	\tailA=(a_{ij})_{i,j=1,\dots,K}\text{ with }a_{ij}=\begin{cases*}
		1,\text{ if }i\leq j\\
		0,\text{ otherwise}
	\end{cases*}\text{,\quad i.e., }
	\tailA=%\begin{pmatrix}
	%	1 & 1 & \cdots & 1\\
	%	0 & \ddots &\ddots & \vdots\\
	%	\vdots & \ddots &\ddots & 1\\
	%	0&\cdots &0&1
	%\end{pmatrix}
	\begin{tikzpicture}[baseline=(current bounding box.center)]
\matrix (m) [matrix of math nodes,nodes in empty cells,right delimiter={)},left delimiter={(},inner sep=-2pt,nodes={inner sep=1ex}]{
1 &  &  & 1\\
0&  & &  \\
&  & & \\
0&  & 0&1 \\
} ;
\draw[loosely dotted,thick] (m-1-1)-- (m-1-4);
\draw[loosely dotted,thick] (m-1-1)-- (m-4-4);
\draw[loosely dotted,thick] (m-1-4)-- (m-4-4);
\draw[loosely dotted,thick] (m-2-1)-- (m-4-3);
\draw[loosely dotted,thick] (m-2-1)-- (m-4-1);
\draw[loosely dotted,thick] (m-4-1)-- (m-4-3);
\end{tikzpicture}.
	\end{equation*} 
%	 is equal to the cone $\tailC=\{y^2-y^1:y^1\leqslant_{\lsum}y^2\}$ induced by the strict partial order $\leqslant_{\lsum}$.
\end{theorem}

\begin{proof}
   First note that $0\in(\tailC\cup\{0\})\cap\hcone(\tailA)$. It thus remains to show that for all $\tilde{u}\in\R^K\setminus\{0\}$, it holds that  $\tilde{u}\in\tailC$ if and only if $\tailA\,\tilde{u}\geqslant 0$.
   
	Now let $\tilde{u}\in\R^K\setminus\{0\}$ with $\tailA\,\tilde{u}\geqslant0$. We define $u\coloneqq0\in\R^K$ and $v\coloneqq\tilde{u}$. Hence, it holds $v-u=\tilde{u}$ and 
	\begin{align*}
		&\tailA\,\tilde{u}\geqslant0\text{ and }\tilde{u}\neq0\\
		\iff\,&\sum_{i=j}^K \tilde{u}_i\geq 0 \quad \text{ for all }j=1,\dots,K,\text{ and }\tilde{u}\neq0\\
		\iff\,&\sum_{i=j}^K v_i\geq \sum_{i=j}^K u_i \quad \text{ for all }j=1,\dots,K,\text{ and }v\neq u\\
		\iff\,& u\leqslant_{\lsum} v\\
        \iff\, & \tilde{u}=v-u\in\tailC.
	\end{align*}
	Thus, we obtain $\hcone(\tailA)\setminus\{0\}= \tailC$, which concludes the proof. %\todo{KK: Das ist doch alles genau dann wenn - dann hat man hier doch "$=$" und ist fertig, oder?}
%
%	For the other direction we assume that  $\tilde{u}\in\tailC$. Thus, there exist $u,v\in\R^K$ such that $v-u=\tilde{u}$ and $u\leqslant_{\lsum}v$. We can conclude
%	\begin{align*}
%		&u\leqslant_{\lsum}v\\
%		\iff\,& \sum_{i=j}^K u_i\leq \sum_{i=j}^K v_i \quad \text{ for all }j=1,\dots,K\\
%		\iff\,& 0 \leq \sum_{i=j}^K v_i -\sum_{i=j}^K u_i\quad\text{ for all }j=1,\dots,K\\
%		\implies\,& 0 \leq \sum_{i=j}^K \tilde{u}_i\quad \text{ for all }j=1,\dots,K\\
%		\iff\,& 0\leqslant \tailA\,\tilde{u}.
%	\end{align*}
%	Consequently, $\tailC\subseteq \hcone(\tailA)\setminus\{0\}$, which concludes the proof.
\end{proof}

%\subsection{Specific Properties of Ordinal Cones}
%------------------------------------------------------------------------------------------------------------------------------

%------------------------------------------------------------------------------------------------------------------------------
%\textcolor{yellow!60!black}{Alternative Beschreibung des resultierenden Kegels über Extremstrahlen (Matrix B)}\\

Theorem~\ref{th:W-M} implies that the closure of the ordinal cone $\tailC\cup\{0\}$, which is a polyhedral cone by Theorem~\ref{thm:A}, must also have a description based on a finite number of extreme rays. Indeed, the following result provides such a description based on exactly $K$ extreme rays.

\begin{theorem}\label{thm:B}
	It holds that $\hcone(\tailA)=\vcone(\tailB)$ for $\tailA$ defined according to Theorem~\ref{thm:A} and $\tailB\in\R^{K\times K}$ given by $\tailB=(b_{ij})_{i,j=1,\dots,K}$ with 
	\begin{equation*}
		b_{ij}=\begin{cases}
			1,&\text{ if }i=j\\
			-1,&\text{ if }i=j-1\\
			0,&\text{ otherwise}
		\end{cases}\text{,\qquad i.e., }\;
	\tailB=%\begin{pmatrix}
%	1 & -1 & 0 & \cdots & 0\\
%	0 & 1 &-1&\ddots & \vdots\\
%	\vdots & \ddots  & \ddots&\ddots & 0\\
%	\vdots & \ddots & \ddots &1 & -1\\
%	0&\cdots&\cdots &0&1
%\end{pmatrix}
\begin{tikzpicture}[baseline=(current bounding box.center)]
\matrix (m) [matrix of math nodes,nodes in empty cells,right delimiter={)},left delimiter={(},inner sep=-2pt,nodes={inner sep=1ex}]{
1 & -1 & 0 &  & 0\\
0 &  & & & \phantom{-1}\\
& \phantom{-1}  & & &0 \\
 &  & & & -1\\
0&& &0&1 \\
} ;
\draw[loosely dotted,thick] (m-1-1)-- (m-5-5);
\draw[loosely dotted,thick] (m-1-2)-- (m-4-5);
\draw[loosely dotted,thick] (m-1-5)-- (m-3-5);
\draw[loosely dotted,thick] (m-1-3)-- (m-3-5);
\draw[loosely dotted,thick] (m-1-3)-- (m-1-5);
\draw[loosely dotted,thick] (m-2-1)-- (m-5-4);
\draw[loosely dotted,thick] (m-2-1)-- (m-5-1);
\draw[loosely dotted,thick] (m-5-1)-- (m-5-4);
%\draw[loosely dotted,thick] (m-5-1)-- (m-5-3);
\end{tikzpicture}.
	\end{equation*}
\end{theorem}
\begin{proof}
	We first show that $\hcone(\tailA)\subseteq\vcone(\tailB)$. Let $d\in\hcone(\tailA)$. Hence, it holds $\tailA\,d\geqq0$ which is equivalent to $\sum_{i=j}^K d_i\geq0$ for all $j=1,\dots,K$. 
	Set $\lambda_j\coloneqq\sum_{i=j}^K d_i\geq0$ and let $B_{j\bullet}$ denote the $j$-th row of $\tailB$, for $j=1,\dots,K$. Then $B_{j \bullet}\,\lambda=\lambda_j-\lambda_{j+1}=d_j$ for $j=1,\dots,K-1$ and $B_{K\bullet}\,\lambda=\lambda_{K}=d_K$. Consequently, we have shown that $d\in\vcone(\tailB)$.
	
	For the other direction, let $d\in\vcone(\tailB)$, i.e.,\ $d=\tailB\,\lambda$ for some $\lambda\geqq0$. The definition of $\tailB$ implies that  $\sum_{i=j}^K B_{i\bullet}\,\lambda=\sum_{i=j}^{K-1}(\lambda_i-\lambda_{i+1})+\lambda_K=\lambda_j$ for all $j=1,\ldots,K-1$ and $B_{K\bullet}\,\lambda=\lambda_K$. Hence, it follows that $\tailA\cdot d=\tailA\cdot (\tailB\,\lambda)=\lambda\geqq0$ and thus $d\in\hcone(\tailA)$, which concludes the proof.
\end{proof}
%\todo{KK: Ab hier wurde die ordinal cone mit $C$ bezeichnet - ich bin dafür, sie wie bisher mit $\tailC$ zu bezeichnen.}
% In the following we will denote $\cone\coloneqq\hcone(A)\setminus\{0\}=\vcone(B)\setminus\{0\}$ as  \emph{ordinal cone}.

%------------------------------------------------------------------------------------------------------------------------------
%\textcolor{yellow!60!black}{Rem: Die Inverse von A ist B (liegt an unserer Skalierung und Sortierung)}\\
Note that these descriptions of the closure of the ordinal cone $\tailC\cup\{0\}$ are not unique. Indeed, both the normal vectors in $\tailA$ as well as the extreme rays in $\tailB$ could be reordered, and they could be multiplied by arbitrary positive scalars without changing the cone that they define. In the particular description given in Theorems~\ref{thm:A} and \ref{thm:B}, however, we observe that the matrix
$\tailB$ is the inverse of the matrix $\tailA$, i.e.,\ $(\tailA)^{-1}=\tailB$. 
%This holds only due to the sorting and scaling we chose for $B$ and because the resulting cone is defined by exactly \kk{$K$} halfspaces or \kk{$K$} extreme rays, respectively. It is possible to define other points on the extreme rays, such that the resulting matrix $\tilde{B}$ is not the inverse of $\kk{\tailA}$.

%------------------------------------------------------------------------------------------------------------------------------
%\textcolor{yellow!60!black}{Bedeutung in unserem Fall: wenn A Matrix der Facetten für Ursprungskegel und B entsprechende Matrix der Extremstrahlen, dann ist $B^\top$ die Matrix der Facetten für den Dualen Kegel und $A^\top$ die Matrix der Extremstrahlen für den Dualen Kegel}\\
\begin{rem}
	It holds $(\hcone(\tailA))^*=\hcone((\tailB)^\top)$ and $(\vcone(\tailB))^*=\vcone((\tailA)^\top)$ since %in our case 
	the normal vectors of the halfspaces given in $\tailA$ are orthogonal to the extreme rays contained in $\tailB$. %\kk{defining $\tailC\cup\{0\}$.}
\end{rem}

\begin{figure}
		\centering
		\subcaptionbox{Pareto Cone 2D \\(Minimization) \label{fig:pareto2D}}
						[.3\linewidth]{\input{Bilder/ParetoCone2D}}
		\subcaptionbox{$\hcone(\tailA)$ 2D\\(Minimization)\label{fig:lastElem2D}}
						[.3\linewidth]{\input{Bilder/LastElemCone2D}}
		\subcaptionbox{$\hcone((\tailA)^\top)$ 2D\\(Maximization)\label{fig:firstElem2D}}
						[.3\linewidth]{\input{Bilder/FirstElemCone2D}}
		\par\bigskip
		\subcaptionbox{Dual Pareto Cone 2D\\(Minimization) \label{fig:paretoDual2D}}
		[.3\linewidth]{\input{Bilder/ParetoCone2D}}
		\subcaptionbox{$(\hcone(\tailA))^*$ 2D\\(Minimization) \label{fig:lastElemDual2D}}
		[.3\linewidth]{\input{Bilder/LastElemDualCone2D}}
		\subcaptionbox{$(\hcone((\tailA)^\top))^*$ 2D\\(Maximization)\label{fig:firstElemDual2D}}
		[.3\linewidth]{\input{Bilder/FirstElemDualCone2D}}
		\par\bigskip
		\subcaptionbox{Pareto Cone 3D\\(Minimization)\label{fig:pareto3D}}
						[.3\linewidth]{\input{Bilder/ParetoCone3D}}
		\subcaptionbox{$\hcone(\tailA)$ 3D\\(Minimization)\label{fig:lastElem3D}}
						[.3\linewidth]{\input{Bilder/LastElemCone3D}}
		\subcaptionbox{$\hcone((\tailA)^\top)$ 3D\\(Maximization)\label{fig:firstElem3D}}
						[.3\linewidth]{\input{Bilder/FirstElemCone3D}}
		\par\bigskip
		\subcaptionbox{Dual Pareto Cone 3D\\(Minimization) \label{fig:paretoDual3D}}
		[.3\linewidth]{\input{Bilder/ParetoCone3D}}
		\subcaptionbox{$(\hcone(\tailA))^*$ 3D\\(Minimization) \label{fig:lastElemDual3D}}
		[.3\linewidth]{\input{Bilder/LastElemDualCone3D}}
		\subcaptionbox{$(\hcone((\tailA)^\top))^*$ 3D\\(Maximization)\label{fig:firstElemDual3D}}
		[.3\linewidth]{\input{Bilder/FirstElemDualCone3D}}
		\caption{Cones and their dual cones}
		\label{fig:cones}
\end{figure}
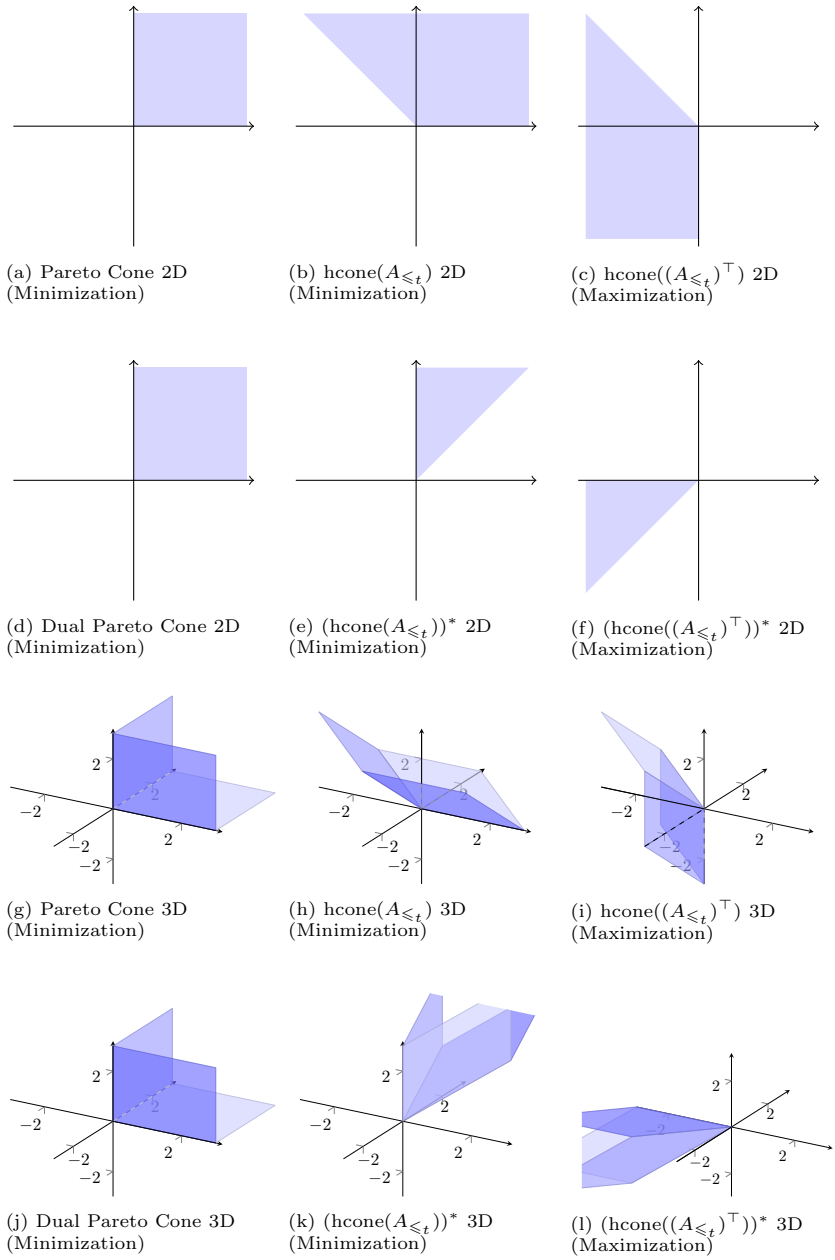

%------------------------------------------------------------------------------------------------------------------------------
%\textcolor{yellow!60!black}{Remark: Ordinal Cone enthält Pareto Cone, Dual Cone ist Teilmenge des Pareto Cone}\\
%Note that the Pareto cone is self dual, i.e., $\pareto^*=\pareto$. 
It is easy to see that the Pareto cone, $\pareto$, is a subset of the ordinal cone, $\tailC$. Moreover, the dual cone of the ordinal cone is a subset of the Pareto cone, i.e.,\ $(\tailC)^*\subseteq \pareto\subseteq\tailC$. This holds since $z\in\pareto$ implies that $z\geqslant 0$ and hence $\tailA\cdot z \geqslant 0$, i.e.,\ $z\in \tailC$. Moreover, $z^*\in(\tailC)^*$ is equivalent to $(z^*)^\top c\geq0$ for all $c\in\tailC$ which implies $z^*_i\geq0$, because the $i$-th unit vector $e_i\in\R^K$ is contained in $\tailC$ for all $i=1,\dots,K$. These cones and their duals are visualized in Figure~\ref{fig:cones}.

%\todo[inline]{KK: Ich würde in Figure~\ref{fig:cones} die beiden Zeilen mit den dualen Kegeln weglassen (subfigures (d),(e),(f),(j),(k),(l)). Oder brauchen wir das später noch? Gegebenenfalls muss die Bildunterschrift noch angepasst werden.}

%------------------------------------------------------------------------------------------------------------------------------
%\textcolor{yellow!60!black}{Bilder 2D und 3D zu den Kegeln (beide Varianten)}\\

\begin{rem}
Note that head-dominance (c.f.\ equation~\eqref{eq:sumFirstElem}), represented by the binary relation $\geqslant_{\fsum}$,  also induces a polyhedral cone. Indeed, we have that $\headC\cup\{0\}=\hcone(\headA)=\vcone(\headB)$, where 
$\headA=\tailA^{\top}$ and $\headB=\tailB^{\top}$.
\end{rem}

%------------------------------------------------------------------------------------------------------------------------------
%Alter Titel: \subsection{Ordinal Cone Optimization Problem and Theoretical Results}
%\subsection{\kk{Ordinal Optimization and Cone Dominance}}\label{subsec:coneOpt}
%\textcolor{yellow!60!black}{Def: Optimalität über Kegel}\\
%The results of Section~\ref{subsec:ordinalcone} above suggest that the ordinal optimization problem \eqref{eq:OCOP} w.r.t.\ tail-dominance (see also Section~\ref{subsec:ocoptail}) can be equivalently formulated as a multi-objective optimization problem w.r.t.\ the ordinal cone $\tailC$. 
Now we can use the ordinal cone $\tailC$ as described in Definition~\ref{def:coneOpt} to reformulate the optimization problem \eqref{eq:OCOP} as follows:

%------------------------------------------------------------------------------------------------------------------------------
%\textcolor{yellow!60!black}{Optimierungsproblem mit Kegel aufstellen}\\

%\kk{Based on Definition~\ref{def:coneOpt} we can restate the ordinal optimization problem w.r.t.\ tail-dominance \eqref{eq:OCOP} as an} \emph{ordinal cone optimization problem} as follows:
%\todo{KK: Auf \eqref{eq:OCOP} verweisen zwei labels: eq:OCOP und eq:COP} %\eqref{eq:COP}
\begin{equation}\label{eq:COP}\tag{OCOP}
    \begin{array}{rl}
         \min\nolimits_{\tailC} \;& c(x)\\
		\text{s.\,t.} \;& x\in \zulTM. \notag
    \end{array}
\end{equation}
Here, $\min_{\tailC}$ denotes the minimization in the sense of Definition \ref{def:coneOpt} for the ordinal cone $\tailC=\hcone(\tailA)\setminus\{0\}$. In other words, the $\tailC$-non-dominated set of problem \eqref{eq:COP} is given by $N(Y,\tailC)$, where $Y=c(X)$. In the following, we use this notation to clearly distinguish between the optimization w.r.t.\ different ordering cones. 

%\todo[inline]{KK: Eigentlich ist hier doch schon alles klar, weil: the ordinal cone $\tailC$ is induced by tail-dominance $\leqslant_{\lsum}$, and conversely, tail-dominance $\leqslant_{\lsum}$ is induced by the ordinal cone $\tailC$. Können wir nicht Def.~20 und Theorem~21 weglassen?}

%------------------------------------------------------------------------------------------------------------------------------
%\textcolor{yellow!60!black}{Transformiertes Problem mit Pareto Kegel aufstellen}\\

\subsection{Bijective Linear Transformation Between Ordinal and Pareto Optimization}

In the previous subsection we showed that tail-dominance, and hence also ordinal dominance due to Lemma~\ref{lem:eff}, can be equivalently described by the ordinal cone $\tailC$.
Moreover, the closure $\tailC\cup\{0\}$ of the ordinal cone is the polyhedral cone $\hcone(\tailA)=\vcone(\tailB)$ that is spanned by $K$ linearly independent extreme rays in $\R^K$, c.f.\ Theorems~\ref{thm:A} and \ref{thm:B}. Since the closure of the Pareto cone $\pareto\cup\{0\}$ is also a polyhedral cone that is  spanned by $K$ linearly independent extreme rays in $\R^K$ (namely the $K$ unit vectors in $\R^K$), there exists a bijective linear transformation that maps the (closure of the) ordinal cone onto the (closure of the) Pareto cone.

We thus define the following \emph{transformed Pareto cone optimization problem} \eqref{eq:POP}
\begin{equation}\label{eq:POP}\tag{TOP}
    \begin{array}{cc}
         \min\nolimits_{\pareto} & \tailA\cdot c(x)\\
		\text{s.\,t.} & x\in \zulTM, \notag
    \end{array}
\end{equation}
where $\min_{\pareto}$ denotes the optimization w.r.t.\ the Pareto cone $\pareto$ according to Definition \ref{def:coneOpt}. Note that the objective vector of problem \eqref{eq:POP} corresponds to the incremental tail counting vector $\tilde{c}(x)= \tailA\cdot c(x)\in\R^K$ introduced in Section~\ref{subsec:problemDef}, that counts in its $j$th component the number of elements of $x$ that are in category $\eta_j$ or worse. %\todo{JS: eigentlich haben wir $\tilde{c}$ genauso schon oben definiert, sollten wir nicht darauf nur verweisen?}
Indeed, for a feasible solution $x=\{e_1,\dots,e_n\}\in X$ we get $\tilde{c}(x)=\sum_{i=1}^n \tilde{c}(e_i)$, where
%We define $\tilde{c}(x)\coloneqq \tailA\cdot c(x)=\sum_{i=1}^n \tilde{c}(e_i)\in\R^K$ for a feasible solution $x=(e_1,\dots,e_n)\in\zulTM$ with 
\begin{equation*}
	\tilde{c}_j(e_i)=\begin{cases}
		1,&\text{ if }\eta_j\preceqq o(e_i)\\
		0,&\text{ otherwise}
	\end{cases}\text{ for all }j=1,\dots,K.
\end{equation*}
Thus, problem~\eqref{eq:POP} is actually a multi-objective optimization problem with $K$ binary objective functions $\tilde{c}_1,\dots,\tilde{c}_K$ defined on the ground set $\dS$, and with feasible set $X\subseteq 2^{\dS}$.
Recall from Section~\ref{subsec:problemDef} that $\tilde{c}_1(e)=1$ for all $e\in\dS$ and hence $\tilde{c}_1(x)$ simply counts the number of elements in a solution $x\in X$. Moreover, the vector $\tilde{c}(e)$ has the consecutive ones property in the sense that whenever a component of $\tilde{c}(e)$ is zero, then all subsequent components of $\tilde{c}(e)$ are also zero. %There are no indices $i>j$ such that $\tilde{c}_i(e)=1$ and $\tilde{c}_j(e)=0$. 
%If we have for example an element $e$ of the best category $\eta_1$ then there have to be zero-entries after the first entry because the element is only considered in the sum over all categories.

As an example, consider the ordinal shortest path problem introduced in  Example~\ref{ex:numRep}. The path $x^1$ consists of the green-dotted edge $e_2$ with $\tilde{c}(e_2)=(1,0,0)^\top$, the orange-dashed edge $e_1$ with $\tilde{c}(e_1)=(1,1,0)^\top$, and the red-solid edge $e_5$ with $\tilde{c}(e_5)=(1,1,1)^\top$. Hence, we compute $\tilde{c}(x^1)=\tilde{c}(e_2)+\tilde{c}(e_1)+\tilde{c}(e_5)=(3,2,1)^\top$, see also Figure~\ref{fig:shortestPath1}.

In order to show that the ordinal counting optimization problem \eqref{eq:COP} (and hence the ordinal optimization problem \eqref{eq:OOP}) can be solved by using the above transformation to the ``standard''  multi-objective optimization problem \eqref{eq:POP}, we use a classical non-dominance mapping result for polyehdral cones. This result can be found in \cite{engau2007domination} and the references therein %\cite{engau2007domination,yu2013multiple,sawaragi1985theory,weidner1990complete,noghin1997relative,cambini2003order,wiecek03}, \todo{KK: Brauchen wir wirklich alle Referenzen?} 
among several others. We include a proof, which is similar to the more general proof in \cite{wiecek03}, for the sake of completeness.

%------------------------------------------------------------------------------------------------------------------------------
%\textcolor{yellow!60!black}{Thm: Transformiertes Problem in Pareto Kegel liefert gleiche Lösungen}\\
%Theorem 3.1.13 in Dissertation.pdf, S. 66, Verweise übernommen aus Dissertation (bei Yu steht in Dissertation 1985????)!
\begin{theorem}[see, e.g., \citealp{engau2007domination}]\label{thm:solution}
	Let $Y\subset\R^K$ be a nonempty set and let $\hcone(A)$ be a cone induced by a matrix $A\in\R^{m\times K}$. Then it holds
	\[
			A\cdot N(Y,\hcone(A)\setminus\{0\})\subseteq N(A\cdot Y,\pareto).
	\]
	If $\rank(A)=K$, then equality holds, i.e., then we have that $A\cdot N(Y,\hcone(A)\setminus\{0\})= N(A\cdot Y,\pareto)$. %holds if the matrix $A$ has maximal rank, i.e., $\rank(A)=K$. 
\end{theorem}
%Beweis anlaog zu dem in \cite{wiecek03}, welcher einen allgemeineren Fall betrachtet, es gibt frühere Quellen!!!
\begin{proof}
	Suppose that $\bar{y}\in Y$ such that $\bar{y}\in N(Y,\hcone(A)\setminus\{0\})$ and $A\cdot \bar{y}\notin N(A\cdot Y,\pareto)$. Then, by Definition~\ref{def:coneOpt}, there exists a $\hat{y}\in Y\setminus\{\bar{y}\}$ such that $A\cdot \hat{y}\in (A\cdot \bar{y}-\pareto)$, i.e.,\ there exists $d\in \pareto$ such that $A\cdot \hat{y}= A\cdot \bar{y}-d$. Hence, it follows that $d=A\cdot \bar{y}-A\cdot \hat{y}=A\cdot(\bar{y}-\hat{y})\geqslant 0$ and thus $\bar{d}\coloneqq\bar{y}-\hat{y}\in\hcone(A)\setminus\{0\}$. %\kk{Moreover,} $\bar{d}\neq0$ since $\bar{y}\neq\hat{y}$. 
	Finally, we can deduce that $\hat{y}\in(\bar{y}-\hcone(A)\setminus\{0\})$, with $\hat{y}\in Y$. But then $\bar{y}\notin N(Y,\hcone(A)\setminus\{0\})$, which contradicts the assumption.
	
	It remains to show that  $A\cdot N(Y,\hcone(A)\setminus\{0\})\supseteq N(A\cdot Y,\pareto)$ if $\rank(A)=K$. Towards this end, suppose that $\bar{y}\in Y$ such that $A\cdot \bar{y}\in N(A\cdot Y,\pareto)$ and $\bar{y}\notin N(Y,\hcone(A)\setminus\{0\})$. Hence, there exists a $d\in\hcone(A)\setminus\{0\}$ such that $\hat{y}=\bar{y}-d\in Y$. This implies $A\cdot \hat{y}=A\cdot\bar{y}-A\cdot d\in A\cdot Y$. From $\rank(A)=K$ and $d\neq 0$ we deduce that $A\cdot d\neq 0$ and thus $A\cdot d\geqslant 0$. Consequently, $(A\cdot\bar{y}-\pareto)\cap A\cdot Y\neq\varnothing$ which contradicts the assumption that $A\cdot \bar{y}\in N(A\cdot Y,\pareto)$.
\end{proof}
	
%------------------------------------------------------------------------------------------------------------------------------
%\textcolor{yellow!60!black}{ Cor: OOP kann mit Pareto Cone gelöst werden!}\\
\begin{theorem}\label{thm:equivalence}
	The set of ordinally efficient solutions for problem \eqref{eq:OOP}, the
	set of tail-efficient (ordinally efficient) solutions of problem \eqref{eq:COP} and the set of Pareto-efficient solutions of problem \eqref{eq:POP} are equal.
\end{theorem}
\begin{proof}
	This follows immediately from Lemma~\ref{lem:eff}, the relation between orders and cones and Theorem~\ref{thm:solution}.
\end{proof}

\begin{figure}
	\centering
	\subcaptionbox{Problem \eqref{eq:COP} \label{fig:3KCOP-shadow}}
	[.4\linewidth]{\input{Bilder/3K-COP-shadow}}
	\subcaptionbox{Problem \eqref{eq:POP} \label{fig:3KPOP-shadow}}
	[.4\linewidth]{\input{Bilder/3K-POP-shadow}}
	\caption{Illustration of the tail-efficient solutions of the instance of \eqref{eq:COP} introduced in Example~\ref{example:cones} (left) and of the respective Pareto-efficient solutions of the  transformed problem \eqref{eq:POP} (right). %Example for the outcome space for an ordinal optimization problem with 3 categories. 
	The dominated areas are shown up to the reference points $(4,4,4)^\top$ (left) and $(5,5,5)^\top$ (right). }
	\label{fig:3K-shadow}
\end{figure}
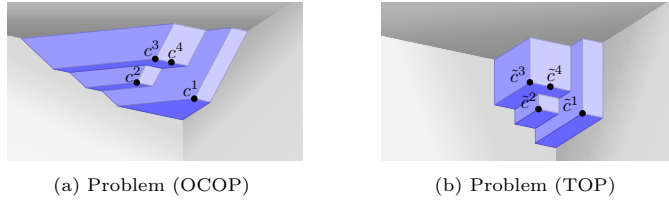

%------------------------------------------------------------------------------------------------------------------------------
%\textcolor{yellow!60!black}{Beispiel mit 3 Kategorien}\\
\begin{example}\label{example:cones}
	Consider an instance of problem \eqref{eq:COP} with $K=3$ categories that has four feasible counting vectors $c^1=(3,1,0)^\top$, $c^2=(0,2,1)^\top$, $c^3=(0,0,2)^\top$ and $c^4=(1,0,2)^\top$. The transformation to problem \eqref{eq:POP} yields the corresponding incremental tail counting vectors as $\tilde{c}^1=(4,1,0)^\top$, $\tilde{c}^2=(3,3,1)^\top$, $\tilde{c}^3=(2,2,2)^\top$ and $\tilde{c}^4=(3,2,2)^\top$. The outcome spaces of both formulations  are depicted in Figure~\ref{fig:3K-shadow} together with the dominance cones $\tailC$ and $\pareto$, respectively.
\end{example}

%%%%%%%%%%%%%%%%%%%%%%%%%%%%%%%%%%%%%%%%%%%%%%%%%%%%%%%%%%%%%%%%%%%%%%%%%%%%%%%%%%%%%%%%%%%%%%%%%%%%%%%%%%%%%%%%%%%%%%%%%%%%%%%%

%%%%%%%%%%%%%%%%%%%%%%%%%%%%%%%%%%%%%%%%
\section{Solution Strategies}\label{sec:solution}
%%%%%%%%%%%%%%%%%%%%%%%%%%%%%%%%%%%%%%%%
%\begin{itemize}
%	\item Algorithmus
%	\item Anmerkung: Alle Eigenschaften der Probleme bleiben erhalten/warum Bellman'sPrinciple angewendet werden kann
%	\item Anmerkung: man findet weniger Lösungen wenn man die Zielfunktion als ordinal betrachtet
%	\item Ordinal Space Decomposition
%\end{itemize}

%------------------------------------------------------------------------------------------------------------------------------
%\textcolor{yellow!60!black}{Algorithmus}\\

In this section, we discuss a generic algorithmic framework for ordinal optimization problems that takes advantage of the close relationship to multi-objective optimization problems. Since the weighted sum scalarization is a popular approach in multi-objective optimization, the interpretation of weights in the context of ordinal optimization and Pareto optimization is analyzed in more detail, and their relation to numerical representations is discussed.
%\todo[inline]{JS: Hier schon sagen, dass die Gewichte der Skalarisierung in Bezug zu den Numerischen Representationen stehen?}

\subsection{Ordinal Optimization by Pareto Transformation}

From the theory above it follows that we can solve the problems \eqref{eq:COP} and \eqref{eq:OOP} by solving the transformed problem \eqref{eq:POP}, which is a standard multi-objective combinatorial problem w.r.t.\ Pareto optimality. After the computation of the Pareto-efficient set of problem \eqref{eq:POP}, or of a minimal complete Pareto-efficient set, respectively, it is necessary to re-compute the corresponding outcome vectors of either problem \eqref{eq:COP} or \eqref{eq:OOP}. In this context, a minimal complete Pareto-efficient set of \eqref{eq:POP} is a subset of the Pareto-efficient set that contains one Pareto-efficient solution for each Pareto-non-dominated outcome vector. We refer to \cite{serafini87some} for different solution concepts in multi-objective optimization.
The efficient sets of the problems \eqref{eq:POP}, \eqref{eq:COP} and \eqref{eq:OOP} are equal and will be denoted by $\effM$ in the following. %This notation is meaningful as the problems have the same efficient set. 
The respective non-dominated sets are
denoted by $\outNdD^{\TOP}\coloneqq N(\tilde{c}(X),\pareto)$, $\outNdD^{\OCOP}\coloneqq N(c(X),\tailC)$ and $\outNdD^{\OOP}$, respectively. A procedure for the computation of the efficient set and the respective non-dominated sets based on this Pareto transformation is outlined in Algorithm~\ref{alg:transformation}.

\begin{algorithm}
	\SetAlgoLined 
	\KwIn{feasible set $\zulTM\subseteq2^\dS$ and ordinal function $o:\dS\to \mathcal{C}$}
	\KwOut{efficient set $\effM$ and non-dominated sets $\outNdD^{\OCOP}$ and $\outNdD^{\OOP}$}
	\caption{Ordinal optimization by Pareto transformation (OOPT)}\label{alg:transformation}
	Compute $c(x)$ for all $x\in\zulTM$ \tcp*{compute counting objective $c$}
	$\effM\coloneqq\min_{\pareto}\{ \tailA\cdot c(x) \colon x\in\zulTM\}$ \tcp*{solve lin.~transf.~\eqref{eq:POP}}
	%\\ $\Rightarrow$ Pareto-efficient solutions $$ and the non-dominated set $\outNdD^{TOP}$\;
	$\outNdD^{\OCOP}\coloneqq c(\effM)$ \tcp*{map efficient set to \ldots} 
	$\outNdD^{\OOP}\coloneqq o(\effM)$\tcp*{\ldots\ resp.\ obj.\ spaces} %for all $x\in\effM$\;
	\Return efficient set $\effM$ and non-dominated sets %$\outNdD^{\TOP}$, 
	$\outNdD^{\OCOP}$ and $\outNdD^{\OOP}$
\end{algorithm}

%------------------------------------------------------------------------------------------------------------------------------
%\textcolor{yellow!60!black}{Anmerkung: Alle Eigenschaften der Probleme bleiben erhalten/warum Bellman'sPrinciple angewendet werden kann}\\
Note that the structural properties of the problems \eqref{eq:COP} and \eqref{eq:OOP} are preserved by the transformation as we do not change the feasible set and as the transformation of the objective function is linear and bijective. In particular, combinatorial solution strategies, like e.g.\ Bellman's principle of optimality for knapsack problems, can be applied in step~2 of Algorithm~\ref{alg:transformation} to efficiently compute $\effM$. Ordinal optimization is thus in general no more complex that standard multi-objective optimization.

%------------------------------------------------------------------------------------------------------------------------------
%\textcolor{yellow!60!black}{Anmerkung: man findet weniger Lösungen wenn man die Zielfunktion als ordinal betrachtet}\\
%For problem \eqref{eq:COP} we could search for the non-dominated vectors in sense of the ordinal cone or we could use the Pareto cone on the same outcome vectors. In the later case, we would receive an equal sized or a larger non-dominated set, because the Pareto cone is a subset of the ordinal cone. 

%-------------------------------------------------------------------------------------------------------------------------------
\subsection{Weighted Sum Scalarization and Ordinal Weight Space Decomposition}

In the following we investigate the interrelation between weighted sum scalarizations for \eqref{eq:POP} and \eqref{eq:COP} and numerical representations for \eqref{eq:OOP}. Thereby we rely on the concept of weight space decompositions, which were introduced by \cite{benson00outcome} for multi-objective linear programming and extended to integer linear problems in \cite{przybylski10recursive}.

The \emph{weighted sum scalarization for \eqref{eq:POP}} is 
	\begin{align}\label{p:wstop}\tag{WSTOP($\lambda$)}
		\min \quad& \sum\limits_{i=1}^K\lambda_i\, \tilde{c}_i(x)\\
		\text{s.\,t.}\quad &x\in \zulM \notag
	\end{align}
with $\lambda_i > 0$ for $i=1,\dots,K$ %$\lambda_i\in(0,1)$ 
and $\sum_{i=1}^K \lambda_i=1$. Analogously, the \emph{weighted sum scalarization for \eqref{eq:COP}} can be formulated as
	\begin{align}\label{p:wscop}\tag{WSOCOP($\mu$)}
		\min \quad& \sum\limits_{i=1}^K \mu_i \, c_i(x)\\
		\text{s.\,t.}\quad&x\in \zulM \notag
	\end{align}
with \(\mu\in (\tailC)_s^*\), where $(\tailC)_s^*$ is the strict dual cone of the ordinal cone $\tailC$, and $\sum_{i=1}^K \mu_i=1$. Recall that the strict dual cone $(\tailC)_s^*$ is the interior of the dual cone $(\tailC)^*$, which is visualized in Figures~\ref{fig:lastElemDual2D} and \ref{fig:lastElemDual3D}.

It is a well-known fact that when considering a multi-objective optimization problem, then optimal solutions of weighted sum scalarizations with weighting vectors $\lambda\in\R^K_{>}$ are always Pareto-efficient  \citep[see, e.g.,][]{Ehrg05}. Such solutions are called \emph{supported} efficient solutions. Thus, problem \eqref{p:wstop} always yields Pareto-efficient solutions for problem \eqref{eq:POP}. Since \eqref{eq:COP} can be interpreted as a multi-objective optimization problems w.r.t.\ the ordering cone \(\tailC\), every optimal solution of the associated weighted sum problem \eqref{p:wscop} with weights in the strict dual cone \((\tailC)^*_s\) of \(\tailC\) is ordinally efficient for \eqref{eq:COP} \citep[see, e.g.,][]{engau2007domination}. 

The supported efficient solutions of problems \eqref{eq:POP} and \eqref{eq:COP} are the same, hence there is a one-to-one correspondence between appropriate weighting vectors $\lambda$ and $\mu$. For a given $\lambda\in\R_{>}^K$ and $x\in\zulTM$ we define $\mu_i=\sum_{j=1}^i\lambda_j$ for all $i=1,\dots,K$. Then it holds that
	$$\sum\limits_{i=1}^K\lambda_i\, \tilde{c}_i(x)=\sum\limits_{i=1}^K\lambda_i\sum\limits_{j=i}^K c_j(x)=\sum\limits_{i=1}^K c_i(x)\sum\limits_{j=1}^i\lambda_j=\sum\limits_{i=1}^K c_i(x)\mu_i,$$
	which shows that problems \eqref{p:wscop} and \eqref{p:wstop} have the same objective functions in this case.
%\begin{theorem}\label{thm:lam_mu}\todo{MS: Enthält dieses Theorem eigentlich noch eine neue Aussage? Reicht nicht das Corollary, der Rest folgt doch aus den Eigenschaften der WS, oder?}
   % Let \replaced[id=MS]{\(\lambda\in\R^K_>\) and \(\mu\in(\tailC)^*_s\).}{$\lambda,\mu\in\R^K_>$. }
%	If $\mu_i=\sum_{j=1}^i\lambda_j$ for all $i=1,\dots,K$, then the sets of optimal solutions of \eqref{p:wscop} and \eqref{p:wstop} are equal. 
%\end{theorem}
%\begin{proof}
   % Let $x\in X$ and let $\mu_i=\sum_{j=1}^i\lambda_j$ for all $i=1,\dots,K$. Then it holds that
%	$$\sum\limits_{i=1}^K\lambda_i\, \tilde{c}_i(x)=\sum\limits_{i=1}^K\lambda_i\sum\limits_{j=i}^K c_j(x)=\sum\limits_{i=1}^K c_i(x)\sum\limits_{j=1}^i\lambda_j=\sum\limits_{i=1}^K c_i(x)\mu_i,$$	which shows that problems \eqref{p:wscop} and \eqref{p:wstop} have the same objective functions and the same feasible set.
%\end{proof}
%\todo[inline]{JS: Aus meiner Sicht kann das Theorem raus. Ich fände es allerdings schön, wenn wir die Umrechnung zwischen $\mu$ und $\lambda$ drin lassen und auch die Umformung von den Zielfunktionen in einander.  Ich würde das aber in den Fließtext schreiben, weil es eigentlich klar ist. Man könnte auch noch einen Satz schreiben, dass die effizeinten Lösungen von den beiden gewichtet Summen Problemen gleich sind, weil die zugrunde liegenden Probleme äquivalent sind.}

Note that $\mu_i=\sum_{j=1}^i\lambda_j$ and $\lambda_i>0$ for all $i=1,\dots,K$ implies that $\mu_i<\mu_j$ for all $i< j$, as required. Conversely, weighting vectors \(\mu\in(\tailC)_s^*\) satisfy $\mu_i<\mu_j$ for all $i<j$ and hence yield associated weighting vectors $\lambda\in\R^K_>$ by setting $\lambda_1\coloneqq\mu_1>0$ and $\lambda_i\coloneqq\mu_i-\mu_{i-1}>0$ for all $i=2,\dots,K$. 
%
%\begin{cor}
%	If \replaced[id=MS]{\(\mu\in(\tailC)^*_s\)}{$\mu$ is an element of the strict dual cone of $\hcone(\tailA)$}, i.e., if $\mu_i<\mu_j$ for all $i< j$, then an optimal solution of problem \eqref{p:wscop} is %a non-dominated 
%	an ordinally efficient solution of problem \eqref{eq:COP}.
%\end{cor}
%
%\begin{proof}
% This result follows from Theorems~\ref{thm:equivalence} and \ref{thm:lam_mu} and the discussion above, and also from a more general result that relates strict dual cones with weighted sum scalarizations, see, e.g.\ \cite{engau2007domination}.
%\end{proof}
%
%\todo[inline]{JS: Das Korollar kann aus meiner Sicht raus, das ist trivial und steht auch schon direkt vor Theorem 24}
%
Note also that while the values of $\mu_i=\sum_{j=1}^i\lambda_j$, $i=1,\dots,K$ (for given $\lambda\in\R_{>}^K$) are in general not normalized to satisfy \(\sum_{i=1}^K \mu_i=1\), such weighting vectors $\mu$ can be easily normalized by setting
\[  
  \mu_i\coloneqq\frac{\sum_{j=1}^{i}\lambda_j}{\sum_{\ell=1}^K\sum_{j=1}^\ell \lambda_j}=\frac{\sum_{j=1}^{i}\lambda_j}{\sum_{j=1}^K(K-j+1) \lambda_j}.
\]
Note that this normalization is  applicable since \((\tailC)_s^*\subset \R^K_>\).
%T\deleted{heorem~\ref{thm:lam_mu} together with t}his normalization of weights imply that 
As a consequence, a weight space decomposition for the multi-objective problem \eqref{eq:POP} can be translated into an associated \emph{ordinal weight space decomposition} for the ordinal counting optimization problem \eqref{eq:COP}. In this context, a weight space decomposition subdivides the space of relevant weighting vectors $\lambda\in\R^K_>$ with $\sum_{i=1}^K\lambda_i=1$ into polyhedral cells such that all weighting vectors from the same cell generate the same efficient solution(s).

\begin{figure}
	\centering
	\subcaptionbox{Weight space decomposition \label{fig:WSD}}
	[.45\linewidth]{\input{Bilder/WeightSpaceDec}}
	\subcaptionbox{Ordinal weight space decomposition\label{fig:OSD}}
	[.45\linewidth]{\input{Bilder/OrdinalSpaceDec}}
	\caption{Weight space decomposition and corresponding ordinal weight space decomposition for the shortest path problem given in Example~\ref{ex:numRep}. The efficient solution $x^2$ corresponds to the light grey triangle, both $x^3$ and $x^4$ correspond to the middle grey triangle and $x^5$ corresponds to the dark grey triangle. The values on dashed lines may not be chosen for $\lambda$ and $\mu$, because for the weight space decomposition we assume that $\lambda\in\R^K_{>}$ and $\sum_{i=1}^{3}\lambda_i=1$, and for the ordinal weight space decomposition we require $0<\mu_1<\mu_2<\mu_3$ and $\sum_{i=1}^3\mu_i=1$.
	%$\mu_1\neq \mu_2$ and $\mu_2\neq\mu_3$.
	}
	\label{fig:SD}
\end{figure}
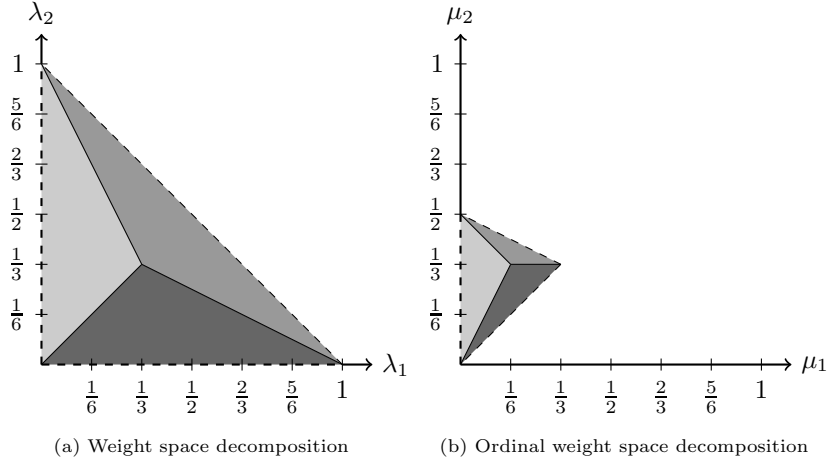

\begin{example}
	In the shortest path problem of Example~\ref{ex:numRep} the solutions $x^2$, $x^3$, $x^4$ and $x^5$ are efficient. In Figure~\ref{fig:SD} the corresponding weight space decomposition is depicted showing the values of $\lambda$ and $\mu$ for which the respective efficient solution is obtained.
\end{example}

From yet another perspective, weighting vectors \(\mu\in(\tailC)_s^*\), i.e., weighting vectors $\mu\in\R^K$ satisfying $0<\mu_i<\mu_j$ for all $i<j$, are related to numerical representations as introduced in Section~\ref{subsec:optConcepts}. Indeed, numerical representations assign a numerical value $\nu(\eta_i)$ to every ordinal category $\eta_i$, $i=1,\dots,K$, such that $\nu(\eta_i)<\nu(\eta_j)$ whenever $i<j$. 
%From $\mu_i=\sum_{j=1}^i\lambda_j$ for all $i=1,\dots,K$ and \(\lambda_j> 0\) for all $j=1,\dots,K$ it follows immediately that $\mu_i<\mu_j$ for $i<j$ and 
Hence we can chose the values \(\mu_i\), $i=1,\dots,K$, equal to the values \(\nu(\eta_i)\) of any numerical representation that satisfies $\nu(\eta_1)>0$. These values can again be normalized without changing the optimal solutions of \eqref{p:wscop} 
%such that its relative size corresponds to the relative size of a numerical representation \(\nu(\eta_i)\), 
by setting
\[
  \mu_i \coloneqq \frac{\nu(\eta_i)}{\sum_{j=1}^K \nu(\eta_j)}, \quad i=1,\dots,K.
\]
It is important to note that this does not imply that numerical representations and weighted sum scalarizations are equivalent. Similarly, it is in general not possible to compute all ordinally efficient solutions of problem~\eqref{eq:OOP} by solving \eqref{p:wscop} for an appropriate $\mu$. To see this, recall that a solution $\xone\in\zulM$ is called ordinally efficient for problem~\eqref{eq:OOP} if and only if there is no $\xtwo\in\zulM$ that ordinally dominates $\xone$, i.e., if  \emph{for every} $\xtwo\in\zulM$ there exists a numerical representation $\vR^{\xtwo}\in\VR$ such that $\vR^{\xtwo}(\xone)\leq\vR^{\xtwo}(\xtwo)$. In contrast, a solution $\xone\in\zulM$ is optimal for problem~\eqref{p:wscop} with appropriate $\mu$ if and only if there exists a numerical representation $\vR^*\in\VR$ such that $\vR^*(\xone)\leq\vR^*(\xtwo)$ \emph{ for all } $\xtwo\in\zulM$.

%-------------------------------------------------------------------------------------------------------------------------------
%\textcolor{yellow!60!black}{Rem: non-supported Lösungen}\\
\begin{rem}
	Ordinal optimization problems may have non-supported efficient solutions. This is illustrated in Example~\ref{ex:nonsupported}. Hence, we can not expect to determine 
	%it is in general not possible to compute 
	all efficient solutions with the weighted sum method.
\end{rem}
	
%	\todo[inline]{JS: Sollten wir aus dem folgenden nicht ein Beispiel machen?}
\begin{example}\label{ex:nonsupported}
	Consider an instance with two categories and  three efficient solutions \(\xone,\xtwo,\xthree\) with counting vectors \(c(\xone)=(3,1)^\top\), \(c(\xtwo)=(5,0)^\top\) and \(c(\xthree)=(0,2)^\top\). The incremental tail counting vectors in the transformed problem \eqref{eq:POP} are $\tilde{c}(\xone)=(4,1)^\top$, $\tilde{c}(\xtwo)=(5,0)^\top$ and $\tilde{c}(\xthree)=(2,2)^\top$, respectively. Obviously, \(\tilde{c}(\xone)\) is non-dominated in \eqref{eq:POP} but unsupported, and thus \(\xone\) is not optimal for \eqref{p:wstop}, irrespective of the choice of $\lambda\in\R^K_>$.
	Similarly, there is no numerical representation such that $\xone$ is simultaneously better than $\xtwo$ \emph{and} $\xthree$, i.e., there is no numerical representation $\nu$ such that
	$\vR(\xone)\leq\vR(\xtwo)$ \emph{and} $\vR(\xone)\leq\vR(\xthree)$. Indeed, the numerical values of the points are $\vR(\xone)=3\, \vR(\eta_1)+\vR(\eta_2)$, $\vR(\xtwo)=5\, \vR(\eta_1)$ and $\vR(\xthree)=2\, \vR(\eta_2)$. 
	$\vR(\xone)\leq\vR(\xtwo)$ implies $\vR(\eta_2)\leq 2\, \vR(\eta_1)$ and $\vR(\xone)\leq\vR(\xthree)$ implies $3\, \vR(\eta_1)\leq\vR(\eta_2)$, which is a contradiction to $\nu(\eta_1)<\nu(\eta_2)$ and % \todo{KK: Sonst wäre $\nu=(0,0)$ eine Lösung} 
	$\vR(\eta_i)\geq 0$ for $i=1,2$.
	However, neither $\xtwo$ nor $\xthree$ ordinally dominate $\xone$, i.e., \emph{neither} $\xtwo$ nor $\xthree$ yield a better objective value for \emph{every} numerical representation.
	
	%As an example for a non-supported solution we consider the following non-dominated points: $\tilde{c}(\bar{x})=(4,1)^\top$, $\tilde{c}(x_1)=(5,0)^\top$ and $\tilde{c}(x_2)=(2,2)^\top$. Hence, it holds $c(\bar{x})=(3,1)^\top$, $c(x_1)=(5,0)^\top$ and $c(x_2)=(0,2)^\top$. We are looking for $\vR$ such that $\vR(\bar{x})\leq\vR(x_i)$ holds for $i=1,2$. The numerical values of the points are $\vR(\bar{x})=3\, \vR(\eta_1)+\vR(\eta_2)$, $\vR(x_1)=5\, \vR(\eta_1)$ and $\vR(x_2)=2\, \vR(\eta_2)$. $\vR(\bar{x})\leq\vR(x_1)$ implies $\vR(\eta_2)\leq2\, \vR(\eta_1)$ and $\vR(\bar{x})\leq\vR(x_2)$ implies $3\, \vR(\eta_1)\leq\vR(\eta_2)$, which is a contradiction for $\vR(\eta_i)\geq 0$ for $i=1,2$.
\end{example}

%%%%%%%%%%%%%%%%%%%%%%%%%%%%%%%%%%%%%%%%%%%%%%%%%%%%%%%%%%%%%%%%%%%%%%
\section{Multi-objective Ordinal  Optimization}\label{sec:extensions}
%%%%%%%%%%%%%%%%%%%%%%%%%%%%%%%%%%%%%%%%%%%%%%%%%%%%%%%%%%%%%%%%%%%%%%%%
%------------------------------------------------------------------------------------------------------------------------------
%\todo[inline]{JS: wenn wir eine ordinale mit einer weiteren Zielfunktion kombinieren wollen, kann man das verbinden. Beispielsweise beim shortest Path Problem könnte man statt in c die Anzahl der Kanten zu zählen, die Länge der Kanten aufsummieren. Frage: Mit nur einer ordinalen Zielfunktion sind matroid Probleme mit dem Greedy Algorithmus lösbar. Nimmt man eine zusätzliche Zielfunktion dazu, dann geht das nicht mehr - Kann man einen effizienten Algorithmus für den Fall formulieren? Müsste man sich nicht auch beim shortest path problem zu Nutze machen können, wenn man nur eine ordinale Zielfunktion betrachtet und keine weitere klassische Zielfunktion?}
%\textcolor{yellow!60!black}{Erweiterung auf den Fall mit mehr Zielfunktionen (Normale und Ordinale)}\\

\subsection{Conflicting Real-valued Objectives and Ordinal Objectives}

The results of Section~\ref{sec:cone} can be extended to multi-objective optimization problems that combine a finite number of $p$ ``standard'' real-valued objective functions
$w^j: \zulTM\to\R$ with a finite number of $r$ ordinal objective functions $o^l:\zulTM\to \mathcal{C}^l$, $l=1,\dots,r$, that are in mutual conflict. The number of categories in the $l$-th ordinal objective function is denoted by $K_l$, i.e.,   $\mathcal{C}^l=\{\eta^l_1,\ldots,\eta^l_{K_l}\}$ for  $l=1,\dots,r$.

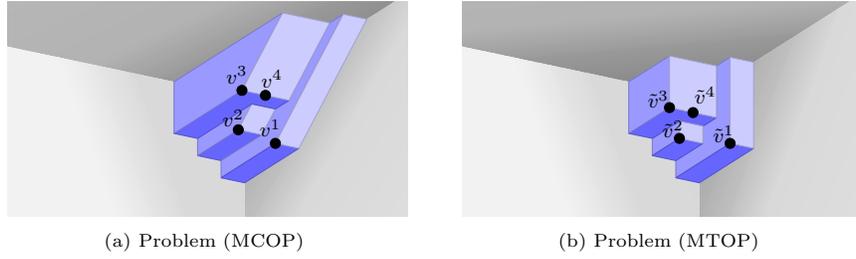
\begin{figure}
	\centering
	\subcaptionbox{Problem \eqref{eq:MCOP} \label{fig:2K+wCOP}}
	[.45\linewidth]{\input{Bilder/2K+w-COP}}
	\hspace{1em}
	\subcaptionbox{Problem \eqref{eq:MPOP} \label{fig:2K+wPOP}}
	[.45\linewidth]{\input{Bilder/2K+w-POP}}
	\caption{Original and transformed outcome space for the multi-objective problem with one real-valued and one ordinal objective function %with two categories 
	introduced in Example~\ref{example:MCOP}. In both figures the upper corner of the bounding box is located at the point $(5,5,5)^\top$.}
	\label{fig:2K+w}
\end{figure}

For a feasible solution $x=\{e_1,\dots,e_n\}\in\zulTM$, we assume that $w^j(x)\coloneqq \sum_{i=1}^n w^j(e_i)$, $j=1,\dots,p$.  
Moreover, $o^l(x)=\sort(o^l(e_1),\dots,o^l(e_n))$ for $l=1,\dots,r$. 
 This leads to the \emph{multi-objective ordinal optimization problem with additional cost functions} \eqref{eq:MOOP}:
\begin{equation}\label{eq:MOOP}\tag{MOOP}
	\begin{array}{ll}
		\min_{\pareto}&(w^1(x),\dots,w^p(x))^{\top}\\
		\min_{\preceqnu} & o^1(x)\\
		\multicolumn{1}{c}{\vdots}&\\
		\min_{\preceqnu} & o^r(x)\\
		\text{s.\,t.} & x\in \zulTM.
	\end{array}
\end{equation}
By replacing the ordered vectors $o^l(x)$ for $l=1,\dots,r$ by the counting vectors $c^l(x)$ for $l=1,\dots,r$
we get a corresponding \emph{multi-objective ordinal counting optimization problem with additional cost functions} \eqref{eq:MCOP}:
\begin{equation}\label{eq:MCOP}\tag{MCOP}
\begin{array}{ll}
    \min\nolimits_{\pareto} & \bigl(w^1(x),\dots,w^p(x)\bigr)^{\top} \\
    \min_{\tailC} & c^1(x)\\
		\multicolumn{1}{c}{\vdots}&\\
    \min_{\tailC} & c^r(x)\\
	%\min\nolimits_\cone \quad& \bigl(c^1(x),\dots,c^r(x)\bigr)^{\kk{\top}} \label{eq:MCOP}\tag{MCOP}\\
	\text{s.\,t.}  & x\in \zulTM.
\end{array}
\end{equation}
We denote the concatenated outcome vectors of \eqref{eq:MCOP} as 
$$v(x) \coloneqq \bigl(w^1(x),\dots,w^p(x),(c^1(x))^{\top},\dots,(c^r(x))^{\top} \bigr)^{\top}\in\R^{p+\tilde{r}},$$ 
where $\tilde{r}\coloneqq\sum_{l=1}^r  K_l$. Then problem \eqref{eq:MCOP} can be transformed into an equivalent standard multi-objective optimization problem w.r.t.\ Pareto dominance using a linear transformation that is defined by the block diagonal matrix % such that we can apply Pareto dominance to all objectives by using matrix
\[
  \tilde{A} \coloneqq \begin{pmatrix} I_{p\times p} & & & \\ & {\tailA^1}&&\\ &&\ddots \\&&& {\tailA^r}\end{pmatrix}.
\]
Here, $I_{p\times p}\in\R^{p\times p}$ denotes the identity matrix and $\tailA^l$ is the transformation matrix corresponding to the objective $c^l$ for $l=1,\dots,r$, c.f.\ Theorem~\ref{thm:A}.
% $\tilde{A}$ which has on its diagonal first a identity matrix $I\in\R^{p\times p}$ and then the matrices $A_l$ corresponding to the objectives $c_l$ for $l=1,\dots,r$ and all other entries are zero. 
Thus, we get the \emph{multi-objective transformed Pareto cone optimization problem} \eqref{eq:MPOP}
\begin{align}
    \min\nolimits_{\pareto} \quad&  \tilde{A}\cdot v(x) \label{eq:MPOP}\tag{MTOP}\\
		\text{s.\,t.} \quad& x\in \zulTM . \notag
\end{align}
Now, problem \eqref{eq:MOOP} or, equivalently, problem \eqref{eq:MCOP} can be solved by using a simple adaptation of Algorithm~\ref{alg:transformation}, c.f.\  Section~\ref{sec:solution}.
%------------------------------------------------------------------------------------------------------------------------------
%\textcolor{yellow!60!black}{Beispiel mit 2 Kategorien und einer normalen Zielfunktion}\\
\begin{example}\label{example:MCOP}
	We consider a problem of type \eqref{eq:MCOP} with one real-valued objective $w$ an one counting objective $c$ with $K=2$ categories (i.e., $p=r=1$). Consider an instance with four feasible outcome vectors $v=(w,c_1,c_2)^{\top}$ given by  $v^1=(4,1,0)^\top$, $v^2=(3,2,1)^\top$, $v^3=(2,0,2)^\top$ and $v^4=(3,0,2)^\top$. Then the corresponding outcome vectors of problem \eqref{eq:MPOP}, $\tilde{v}^i=\tilde{A}\,v^i$ for $i=1,\dots,4$,  
	are obtained as $\tilde{v}^1=(4,1,0)^\top$, $\tilde{v}^2=(3,3,1)^\top$, $\tilde{v}^3=(2,2,2)^\top$ and $\tilde{v}^4=(3,2,2)^\top$. In this case, the transformation matrix is given by
	$$ \tilde{A}=\begin{pmatrix}
	1& 0& 0\\
	0& 1 &1 \\
	0& 0& 1
	\end{pmatrix}.$$ The feasible points and the dominated volumes in the respective outcome spaces are depicted in Figure~\ref{fig:2K+w} for both problems, \eqref{eq:MCOP} and \eqref{eq:MPOP}.
\end{example}

\subsection{Coherent Real-valued Objectives and Ordinal Objectives}

In some practical applications, the elements of $\dS$ have a real-valued cost (e.g., the length of an edge) and an associated category (e.g., the safety of the corresponding road segment for a cyclist) such that the real-valued cost is, rather than in conflict, coherent with the respective category. This situation is illustrated at the following example:   
%    There are optimization models where the real-valued objective is not considered as an independent optimization goal but rather as an attribute compatible with the ordinal objective value. To clarify this point
	%Depending on the context, the method described above may not be meaningful. For an example 
\begin{example}\label{ex:sec}
	Consider the shortest path problem shown in Figure~\ref{fig:extension}. Let $w(e)$ denote the length of an edge $e$ and let $o(e)$ denote its safety: dotted green edges are save and in category $\eta_1$, while solid red edges are insecure and in category $\eta_2$. Then, irrespective of the number of edges contained in the respective paths, the path $x^1=\{e_1,e_2,e_3\}$ should be preferred over the path $x^2=\{e_4,e_5,e_6\}$ since the total weights are equal $w(x^1)=w(x^2)=10)$, and the red sub-path in $x^1$ has a smaller weight than that of $x^2$. 
	In this sense, the weight or length of an edge can be interpreted as an attribute of its respective category.
% 	Nevertheless, in the sense of 
	However, \(x^1\) is dominated by \(x^2\) w.r.t.\ problem~\eqref{eq:MOOP}. 
\end{example}

%\todo{Verweis Fahrrad?}
To model the situation where a real-valued objective function $w:S\to \R$ is in accordance with an ordinal objective function $o:S\to \mathcal{C}$ with $K$ categories, i.e., the situation where the weight $w(e)$ reflects the multiplicity with which the category $o(e)$ of the element $e$ is to be counted, we introduce a \emph{weighted counting vector}
%\emph{category conform weight function}  $c_w:S\to\R^K$ with
	\begin{align*}
		c^w_i(e)\coloneqq\begin{cases}
			w(e) &  \text{if~} o(e)=\eta_i\\
			0& \text{otherwise.}
		\end{cases}
	\end{align*}%\todo{MS: Notation \(c_{w,i}\) oder \(c_i^w\) oder ganz anders?}
% Assume for example that for an element $e\in S$ it holds $w(e)=7$, $o(e)=\eta_2$ and $K=4$.
The basic idea of this concept is also used in the risk-aware bicycle routing application \citetalias{geovelo}, which takes, besides the route length, also the total length of unsafe route segments into account.
For example, when \(K=4\),  \(w(e)=7\) and \(o(e)=\eta_2\), then $c^w(e)=(0,7,0,0)^\top$. 
%Taking the sum over all objective function values of the elements of a feasible solution $x=\{e_1,\dots,e_n\}\in\zulTM$ leads to $c_w(x)=\sum_{i=1}^n c_w(e_i)$. 
The weighted counting objective of a feasible solution \(x=\{e_1,\dots,e_n\}\in\zulTM\) equals the sum of the weighted counting vectors %category conform weight function values 
of all elements in \(x\), i.e., $c^w(x)=\sum_{i=1}^n c^w(e_i)$. Thereby, the
$i$-th component of $c^w(x)$ corresponds to the total weight of the elements in $x$ that are in category $\eta_i$, $i=1,\dots,K$. Now the weighted counting vector can be handled analogously to the counting vector. Indeed, as in the previous chapter we consider the transformation $\tilde{c}^w_i(x)\coloneqq\sum_{j=i}^K c^w_j(x)=\tailA\cdot c^{w}(x)$ to obtain the \emph{weighted transformed Pareto cone optimization problem} %with weighted elements}
\begin{equation}\label{eq:COPWE}\tag{WTOP}%\tag{COPWE}
	\begin{array}{rl}
		\min_{\pareto}&\tilde{c}^w(x)\\
		\text{s.\,t.} & x\in \zulTM
	\end{array}
\end{equation}
w.r.t.\ the concept of Pareto optimality, that can be solved with the methods developed in the preceding sections.

\subsection{Modelling Aspects}

We emphasize that it depends on the context of the respective application whether the multi-objective model~\eqref{eq:MPOP} or the aggregated model \eqref{eq:COPWE} is more suitable. The following example illustrates that the aggregated model \eqref{eq:COPWE} is meaningful whenever $w$ and $c$ are interrelated and coherent objectives, while
the multi-objective model \eqref{eq:MPOP} is particularly useful for unrelated or incompatible objectives.

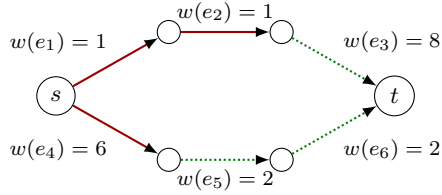
\begin{figure}
    \centering\small
		\include{Bilder/graphExtension}
%  		\vspace*{0.5cm}
		\begin{tabular}{c|cccc}
		    \toprule
			&$w$&$o$&$c^w$&$\tilde{c}^w$\\
			\midrule
			\rule{0pt}{25pt}
			$x^1=\{e_1,e_2,e_3\}$&$10$&$\begin{pmatrix} \eta_1\\ \eta_2\\ \eta_2 \end{pmatrix}$&$\begin{pmatrix} 8\\ 2 \end{pmatrix}$&$\begin{pmatrix} 10\\ 2 \end{pmatrix}$\\
			\rule{0pt}{25pt}
			$x^2=\{e_4,e_5,e_6\}$&$10$&$\begin{pmatrix} \eta_1\\ \eta_1\\ \eta_2 \end{pmatrix}$&$\begin{pmatrix} 4\\ 6 \end{pmatrix}$&$\begin{pmatrix} 10\\ 6 \end{pmatrix}$\\\bottomrule
		\end{tabular}
	\caption{Instance of a shortest path problem. A dotted-green edge is in the best category $\eta_1$ and the solid-red edges are in the worst category $\eta_2$. The possible $s$-$t$-paths $x^1$ and $x^2$ and their different objective function vectors are given. }
	\label{fig:extension}
\end{figure}

\begin{example}
	Consider again the shortest path problem depicted in Figure~\ref{fig:extension}. Obviously, the path $x^1=\{e_1,e_2,e_3\}$ is the unique efficient solution for problem~\eqref{eq:COPWE}, while the path $x^2=\{e_4,e_5,e_6\}$ is the unique efficient solution for problem~\eqref{eq:MOOP}. 
		
	%	\js{If $w(e)$ represents the length of the road $e$, then it makes sense, that the upper path $x^1$ is the better one, because the total length of red edges is shorter in path $x^1$ than in path $x^2$. Nevertheless, the number of red edges is larger in $x^1$ than in $x^2$ and the total length is the same for both paths.}
	%
	%\js{This may be meaningful if} the categories represent the quality of the road, i.e., green edges correspond to asphalt roads in a good condition, while red edges are roads with many potholes.
    Whether $x^1$ or $x^2$ are actually preferred thus depends on the interpretation of the weights and of the ordinal categories. Towards this end, suppose that,
    as in Example~\ref{ex:sec}, the category of an edge corresponds to its security level. 
    
    First, consider the case that the real-valued objective \(w(e)\) represents the length of the edge \(e\) as in Example~\ref{ex:sec}. Then both objectives are interrelated and hence model \eqref{eq:COPWE} is appropriate. 
    
    If, on the other hand, $w(e)$ represents the toll of the edge or road $e$, then this real-valued objective is not an attribute of the corresponding category. In other words, both objectives are potentially conflicting and not coherent. Hence, in this case the path  $x^2$ is preferred since the total amount of toll is the same for both paths, but the second path has more green and fewer red edges.
\end{example}

%%%%%%%%%%%%%%%%%%%%%%%%%%%%%%%%%%%%%%%%%%%%%%%%%%%%%%%%%%%%%%%%%%%%%%%%%%%%%%%%%%%%%%%%%%%%%%%%%%%%%%%%%%%%%%%%%%%%%%%%%%%%%%%

%%%%%%%%%%%%%%%%%%%%%%%%%%%%%%%%%%%%%%%%
\section{Conclusion}\label{sec:conclusion}
%%%%%%%%%%%%%%%%%%%%%%%%%%%%%%%%%%%%%%%%
In this paper we investigate ordinal combinatorial optimization problems. We describe  different optimality concepts for ordinal objective functions, namely ordinal optimality, which was introduced first by \cite{SCHAFER:knapsack}, as well as tail- and head-optimality. We prove that all three concepts are equivalent if all feasible solutions have the same length. In general, only ordinal optimality and tail-optimality are equivalent.

We provide alternative descriptions of these three optimality concepts based on associated ordering cones. Using the fact that  ordinal optimality and tail-optimality are equivalent, and that  tail-optimality can be represented by a polyhedral cone with $K$ extreme rays in $\R^K$, we show that ordinal optimization problems can be transformed into equivalent multi-objective optimization problems with binary cost coefficients. The transformation is realized by a bijective linear mapping. The resulting problem can be solved with standard methods from multi-objective optimization, and hence ordinal optimization is as easy or hard as the associated, ``standard'' multi-objective problems. For example, ordinal knapsack problems and ordinal shortest path problems can be solved by multi-objective dynamic programming, using Bellman's principle of optimality. %Additionally, we explain the relation between numerical representations and weight space decompositions of the multi-objective optimization problem.

The results can be extended to problems with more than one objective function. We suggest two modelling approaches to combine an ordinal objective function with a real-valued objective function. While in the first approach all objectives are considered in a standard multi-objective setting, the second approach allows to model interrelated and coherent objective functions, where the real-valued objective is interpreted as an attribute of the respective category in the ordinal objective. 
%Which one is preferred depends on the context and whether the objective functions are interrelated or not. 

Future work should focus on the development of tailored optimization algorithms for the associated multi-objective optimization problems that exploit the fact that these problems have binary cost coefficients. Moreover, specific combinatorial problems like, for example, shortest path, knapsack, assignment and general routing and network flow problems should be analyzed both with coherent and with conflicting ordinal and real-valued objective functions.
%\begin{itemize}
%    \item lässt sich die spezielle Struktur bei bestimmten kombinatorischen Problemen (kürzeste Wege, Matroid Probleme mit coherent Real-valued and ordinal Objectives) gezielt nutzen?
%    \item 
%\end{itemize}

\section*{Acknowledgment}
The authors thankfully acknowledge financial support by Deutsche Forschungsgemeinschaft, project number~KL~1076/11-1.

\bibliography{literaturMatroid}
\end{document}

%% file: Bilder/bspDirGraph.tex
%\documentclass {article}
%\usepackage{tikz}
%\usetikzlibrary{graphs,quotes}

%\begin{document}

\begin{tikzpicture}[scale=1.5,every edge quotes/.append style={font=\footnotesize}]%[shorten >=1pt,->,draw=black!75, node distance=\layersep,scale=0.28]
\node[draw,circle] (1) at (0,0)[] {$s$};
\node[draw,circle] (2) at (1,1) {};
\node[draw,circle] (3) at (1,0) {};
\node[draw,circle] (4) at (2,0) {$t$};
\node[draw,circle] (5) at (1,-1) {};

%\node (BotLeft) at (-2,-6) {};
%\node (BotLeft2) at (-1.75,-4.75) {};

%%graph G
\graph {
	(1) ->["$e_1$",-latex,dashed,thick,draw=yellow!40!red] (2);
	(2) ->["$e_2$",-latex,densely dotted,thick,swap,draw=black!50!green] (3);
	(2) ->["$e_3$",-latex,dashed,thick,draw=yellow!40!red] (4);
	(1) ->["$e_4$",-latex,densely dotted,thick,swap,draw=black!50!green] (3);
	(3) ->["$e_5$",-latex,thick,swap,draw=black!40!red] (4);
	(1) ->["$e_6$",-latex,dashed,swap,thick,draw=yellow!40!red] (5);
	(3) ->["$e_7$",-latex,densely dotted,thick,swap,draw=black!50!green] (5);
	(5) ->["$e_8$",-latex,dashed,swap,thick,draw=yellow!40!red] (4);
};
\end{tikzpicture}

%\end{document}

%% file: Bilder/numerRepraesentation1-neu.tex
\begin{tikzpicture}[scale=0.4]
	%\draw[very thin,color=gray] (0,8) grid (29,24);
	\node (A) at (1,-0.5) [circle, font=\footnotesize] {$\eta_1$};
	\node (B) at (2,-0.5) [circle, font=\footnotesize] {$\eta_2$};
	\node (C) at (3,-0.5) [circle, font=\footnotesize] {$\eta_2$};
	\node (D) at (4,-0.5) [circle, font=\footnotesize] {$\eta_3$};
	\node (E) at (5,-0.5) [circle, font=\footnotesize] {$\eta_4$};
	\node (F) at (6,-0.5) [circle, font=\footnotesize] {$\eta_4$};
	\node (G) at (7,-0.5) [circle, font=\footnotesize] {$\eta_4$};

	\draw[fill=white, draw =black] (1,1) circle (2ex);
	\draw[fill=white, draw =black] (1,2) circle (2ex);
	
	\draw[fill=black!12, draw =black] (2,1) circle (2ex);
	\draw[fill=black!12, draw =black] (2,2) circle (2ex);
	\draw[fill=black!12, draw =black] (2,3) circle (2ex);
	\draw[fill=black!12, draw =black] (2,4) circle (2ex);
	
	\draw[fill=black!25, draw =black] (3,1) circle (2ex);
	\draw[fill=black!25, draw =black] (3,2) circle (2ex);
	\draw[fill=black!25, draw =black] (3,3) circle (2ex);
	\draw[fill=black!25, draw =black] (3,4) circle (2ex);
	
	\draw[fill=black!43, draw =black] (4,1) circle (2ex);
	\draw[fill=black!43, draw =black] (4,2) circle (2ex);
	\draw[fill=black!43, draw =black] (4,3) circle (2ex);
	\draw[fill=black!43, draw =black] (4,4) circle (2ex);
	\draw[fill=black!43, draw =black] (4,5) circle (2ex);
	\draw[fill=black!43, draw =black] (4,6) circle (2ex);
	\draw[fill=black!43, draw =black] (4,7) circle (2ex);
	
	\draw[fill=black!60, draw =black] (5,1) circle (2ex);
	\draw[fill=black!60, draw =black] (5,2) circle (2ex);
	\draw[fill=black!60, draw =black] (5,3) circle (2ex);
	\draw[fill=black!60, draw =black] (5,4) circle (2ex);
	\draw[fill=black!60, draw =black] (5,5) circle (2ex);
	\draw[fill=black!60, draw =black] (5,6) circle (2ex);
	\draw[fill=black!60, draw =black] (5,7) circle (2ex);
	\draw[fill=black!60, draw =black] (5,8) circle (2ex);
	
	\draw[fill=black!77, draw =black] (6,1) circle (2ex);
	\draw[fill=black!77, draw =black] (6,2) circle (2ex);
	\draw[fill=black!77, draw =black] (6,3) circle (2ex);
	\draw[fill=black!77, draw =black] (6,4) circle (2ex);
	\draw[fill=black!77, draw =black] (6,5) circle (2ex);
	\draw[fill=black!77, draw =black] (6,6) circle (2ex);
	\draw[fill=black!77, draw =black] (6,7) circle (2ex);
	\draw[fill=black!77, draw =black] (6,8) circle (2ex);
	
	\draw[fill=black, draw =black] (7,1) circle (2ex);
	\draw[fill=black, draw =black] (7,2) circle (2ex);
	\draw[fill=black, draw =black] (7,3) circle (2ex);
	\draw[fill=black, draw =black] (7,4) circle (2ex);
	\draw[fill=black, draw =black] (7,5) circle (2ex);
	\draw[fill=black, draw =black] (7,6) circle (2ex);
	\draw[fill=black, draw =black] (7,7) circle (2ex);
	\draw[fill=black, draw =black] (7,8) circle (2ex);
	
	% Achse zeichnen
	\draw[->,thick] (0,0) -- (0,9) node[above, font=\footnotesize] {$\vR(\eta_i)$};
	% Achse beschriften
	\foreach \y in {1,2,3,4,5,6,7,8}\draw (-0.1,\y) -- (0.1,\y) node[left=5pt] {$\scriptstyle\y$};
\end{tikzpicture}

%% file: Bilder/numerRepraesentation2.tex
\begin{tikzpicture}[scale=0.4]
	%\draw[very thin,color=gray] (0,8) grid (29,24);
	\node (A) at (1,-0.5) [circle, font=\footnotesize] {$\eta_1$};
	\node (B) at (2,-0.5) [circle, font=\footnotesize] {$\eta_2$};
	\node (C) at (3,-0.5) [circle, font=\footnotesize] {$\eta_2$};
	\node (D) at (4,-0.5) [circle, font=\footnotesize] {$\eta_3$};
	\node (E) at (5,-0.5) [circle, font=\footnotesize] {$\eta_4$};
	\node (F) at (6,-0.5) [circle, font=\footnotesize] {$\eta_4$};
	\node (G) at (7,-0.5) [circle, font=\footnotesize] {$\eta_4$};

	\draw[fill=white, draw =black] (1,1) circle (2ex);
	\draw[fill=white, draw =black] (1,2) circle (2ex);
	
	\draw[fill=black!25, draw =black] (2,1) circle (2ex);
	\draw[fill=black!25, draw =black] (2,2) circle (2ex);
	\draw[fill=black!25, draw =black] (2,3) circle (2ex);
	\draw[fill=black!25, draw =black] (2,4) circle (2ex);
	
	\draw[fill=black!25, draw =black] (3,1) circle (2ex);
	\draw[fill=black!25, draw =black] (3,2) circle (2ex);
	\draw[fill=black!25, draw =black] (3,3) circle (2ex);
	\draw[fill=black!25, draw =black] (3,4) circle (2ex);
	
	\draw[fill=black!50, draw =black] (4,1) circle (2ex);
	\draw[fill=black!50, draw =black] (4,2) circle (2ex);
	\draw[fill=black!50, draw =black] (4,3) circle (2ex);
	\draw[fill=black!50, draw =black] (4,4) circle (2ex);
	\draw[fill=black!50, draw =black] (4,5) circle (2ex);
	\draw[fill=black!50, draw =black] (4,6) circle (2ex);
	\draw[fill=black!50, draw =black] (4,7) circle (2ex);
	
	\draw[fill=black, draw =black] (5,1) circle (2ex);
	\draw[fill=black, draw =black] (5,2) circle (2ex);
	\draw[fill=black, draw =black] (5,3) circle (2ex);
	\draw[fill=black, draw =black] (5,4) circle (2ex);
	\draw[fill=black, draw =black] (5,5) circle (2ex);
	\draw[fill=black, draw =black] (5,6) circle (2ex);
	\draw[fill=black, draw =black] (5,7) circle (2ex);
	\draw[fill=black, draw =black] (5,8) circle (2ex);
	
	\draw[fill=black, draw =black] (6,1) circle (2ex);
	\draw[fill=black, draw =black] (6,2) circle (2ex);
	\draw[fill=black, draw =black] (6,3) circle (2ex);
	\draw[fill=black, draw =black] (6,4) circle (2ex);
	\draw[fill=black, draw =black] (6,5) circle (2ex);
	\draw[fill=black, draw =black] (6,6) circle (2ex);
	\draw[fill=black, draw =black] (6,7) circle (2ex);
	\draw[fill=black, draw =black] (6,8) circle (2ex);
	
	\draw[fill=black, draw =black] (7,1) circle (2ex);
	\draw[fill=black, draw =black] (7,2) circle (2ex);
	\draw[fill=black, draw =black] (7,3) circle (2ex);
	\draw[fill=black, draw =black] (7,4) circle (2ex);
	\draw[fill=black, draw =black] (7,5) circle (2ex);
	\draw[fill=black, draw =black] (7,6) circle (2ex);
	\draw[fill=black, draw =black] (7,7) circle (2ex);
	\draw[fill=black, draw =black] (7,8) circle (2ex);
	
	% Achse zeichnen
	\draw[->,thick] (0,0) -- (0,9) node[above, font=\footnotesize] {$\vR(\eta_i)$};
	% Achse beschriften
	\foreach \y in {1,2,3,4,5,6,7,8}\draw (-0.1,\y) -- (0.1,\y) node[left=5pt] {$\scriptstyle\y$};
\end{tikzpicture}

%% file: Bilder/numerRepraesentation3.tex
\begin{tikzpicture}[scale=0.4]
	%\draw[very thin,color=gray] (0,8) grid (29,24);
	\node (A) at (1,-0.5) [circle, font=\footnotesize] {$\eta_1$};
	\node (B) at (2,-0.5) [circle, font=\footnotesize] {$\eta_2$};
	\node (C) at (3,-0.5) [circle, font=\footnotesize] {$\eta_2$};
	\node (D) at (4,-0.5) [circle, font=\footnotesize] {$\eta_3$};
	\node (E) at (5,-0.5) [circle, font=\footnotesize] {$\eta_4$};
	\node (F) at (6,-0.5) [circle, font=\footnotesize] {$\eta_4$};
	\node (G) at (7,-0.5) [circle, font=\footnotesize] {$\eta_4$};

	\draw[fill=white, draw =black] (1,1) circle (2ex);
	\draw[fill=white, draw =black] (1,2) circle (2ex);
	
	\draw[fill=white, draw =black] (2,1) circle (2ex);
	\draw[fill=white, draw =black] (2,2) circle (2ex);
	\draw[fill=black!25, draw =black] (2,3) circle (2ex);
	\draw[fill=black!25, draw =black] (2,4) circle (2ex);
	
	\draw[fill=white, draw =black] (3,1) circle (2ex);
	\draw[fill=white, draw =black] (3,2) circle (2ex);
	\draw[fill=black!25, draw =black] (3,3) circle (2ex);
	\draw[fill=black!25, draw =black] (3,4) circle (2ex);
	
	\draw[fill=white, draw =black] (4,1) circle (2ex);
	\draw[fill=white, draw =black] (4,2) circle (2ex);
	\draw[fill=black!25, draw =black] (4,3) circle (2ex);
	\draw[fill=black!25, draw =black] (4,4) circle (2ex);
	\draw[fill=black!50, draw =black] (4,5) circle (2ex);
	\draw[fill=black!50, draw =black] (4,6) circle (2ex);
	\draw[fill=black!50, draw =black] (4,7) circle (2ex);
	
	\draw[fill=white, draw =black] (5,1) circle (2ex);
	\draw[fill=white, draw =black] (5,2) circle (2ex);
	\draw[fill=black!25, draw =black] (5,3) circle (2ex);
	\draw[fill=black!25, draw =black] (5,4) circle (2ex);
	\draw[fill=black!50, draw =black] (5,5) circle (2ex);
	\draw[fill=black!50, draw =black] (5,6) circle (2ex);
	\draw[fill=black!50, draw =black] (5,7) circle (2ex);
	\draw[fill=black, draw =black] (5,8) circle (2ex);
	
	\draw[fill=white, draw =black] (6,1) circle (2ex);
	\draw[fill=white, draw =black] (6,2) circle (2ex);
	\draw[fill=black!25, draw =black] (6,3) circle (2ex);
	\draw[fill=black!25, draw =black] (6,4) circle (2ex);
	\draw[fill=black!50, draw =black] (6,5) circle (2ex);
	\draw[fill=black!50, draw =black] (6,6) circle (2ex);
	\draw[fill=black!50, draw =black] (6,7) circle (2ex);
	\draw[fill=black, draw =black] (6,8) circle (2ex);
	
	\draw[fill=white, draw =black] (7,1) circle (2ex);
	\draw[fill=white, draw =black] (7,2) circle (2ex);
	\draw[fill=black!25, draw =black] (7,3) circle (2ex);
	\draw[fill=black!25, draw =black] (7,4) circle (2ex);
	\draw[fill=black!50, draw =black] (7,5) circle (2ex);
	\draw[fill=black!50, draw =black] (7,6) circle (2ex);
	\draw[fill=black!50, draw =black] (7,7) circle (2ex);
	\draw[fill=black, draw =black] (7,8) circle (2ex);
	
	% Achse zeichnen
	\draw[->,thick] (0,0) -- (0,9) node[above, font=\footnotesize] {$\vR(\eta_i)$};
	% Achse beschriften
	\foreach \y in {1,2,3,4,5,6,7,8}\draw (-0.1,\y) -- (0.1,\y) node[left=5pt] {$\scriptstyle\y$};
\end{tikzpicture}

%% file: Bilder/ParetoCone2D.tex
\begin{tikzpicture}[scale=0.5]
	% cone
	\fill[fill=blue!20,fill opacity=0.8] (0,3) -- (3,3) -- (3,0) -- (0,0)  -- cycle;
	
	% axes    
	\draw[->] (-3.2,0) -- (3.2,0) node[right] {};
	\draw[->] (0,-3.2) -- (0,3.2) node[above] {};
\end{tikzpicture}

%% file: Bilder/LastElemCone2D.tex
\begin{tikzpicture}[scale=0.5]
	% cone
	\fill[fill=blue!20,fill opacity=0.8] (3,0) -- (3,3) -- (-3,3) -- (0,0)  -- cycle;
	
	% axes  
	\draw[->] (-3.2,0) -- (3.2,0) node[right] {};
	\draw[->] (0,-3.2) -- (0,3.2) node[above] {};
\end{tikzpicture}

%% file: Bilder/FirstElemCone2D.tex
\begin{tikzpicture}[scale=0.5]
	% cone
	\fill[fill=blue!20,fill opacity=0.8] (-3,3) -- (-3,-3) -- (0,-3) -- (0,0)  -- cycle;
	
	% axes
	\draw[->] (-3.2,0) -- (3.2,0) node[right] {};
	\draw[->] (0,-3.2) -- (0,3.2) node[above] {};
\end{tikzpicture}

%% file: Bilder/LastElemDualCone2D.tex
\begin{tikzpicture}[scale=0.5]
	% cone
	\fill[fill=blue!20,fill opacity=0.8] (3,3) -- (0,3) -- (0,0)  -- cycle;
	
	% axes   
	\draw[->] (-3.2,0) -- (3.2,0) node[right] {};
	\draw[->] (0,-3.2) -- (0,3.2) node[above] {};

\end{tikzpicture}

%% file: Bilder/FirstElemDualCone2D.tex
\begin{tikzpicture}[scale=0.5]
	% cone	
	\fill[fill=blue!20,fill opacity=0.8] (-3,0) -- (-3,-3) -- (0,0)  -- cycle;
	
	% axes   
	\draw[->] (-3.2,0) -- (3.2,0) node[right] {};
	\draw[->] (0,-3.2) -- (0,3.2) node[above] {};
\end{tikzpicture}

%% file: Bilder/ParetoCone3D.tex
\begin{tikzpicture}[scale=0.65]
	\begin{axis}[
		view/h=30,
		view/v=30,
		xmin=-3, xmax=3.2,
		ymin=-3, ymax=3.2,
		zmin=-3, zmax=3.2,
		colormap={blau}{color=(blue!20) color=(blue!40) color=(blue!60)}
		]
		\addplot3[
		opacity=0.6,
		table/row sep=\\,
		patch,
		patch type=polygon,
		vertex count=4,
		patch table with point meta={%
			% pt1 pt2 pt3 pt4 pt5 cdata
			0 1 4 2 0\\
			0 3 5 2 1\\
			0 1 6 3 2\\
		}]
		table {
			x y z\\
			0 0 0\\%0
			3 0 0\\%1
			0 3 0\\%2
			0 0 3\\%3
			3 3 0\\%4
			0 3 3\\%5
			3 0 3\\%6
		};
		\draw[dashed, black!20] (axis cs: 0,0,0) -- (axis cs: 0,3,0);
		\draw[-, black] (axis cs: 0,0,0) -- (axis cs: 0,0,3);
		%\draw[dashed, gray] (axis cs: 0,0,0) -- (axis cs: 1,0,0);
		\draw[-, black] (axis cs: 0,0,0) -- (axis cs: 3,0,0);
		%\draw[dashed, gray] (axis cs: 0,0,0) -- (axis cs: 0,1,0);
		%\draw[-, black] (axis cs: 0,1,0) -- (axis cs: 0,3.5,0);
	\end{axis}
\end{tikzpicture}

%% file: Bilder/LastElemCone3D.tex
\begin{tikzpicture}[scale=0.65]
	\begin{axis}[
		view/h=30,
		view/v=30,
		xmin=-3, xmax=3.2,
		ymin=-3, ymax=3.2,
		zmin=-3, zmax=3.2,
		colormap={blau}{color=(blue!20) color=(blue!40) color=(blue!60)}
		]
		\addplot3[
		opacity=0.6,
		table/row sep=\\,
		patch,
		patch type=polygon,
		vertex count=4,
		patch table with point meta={%
			% pt1 pt2 pt3 pt4 pt5 cdata
			0 1 4 2 0\\
			0 3 5 2 1\\
			0 1 6 3 2\\
		}]
		table {
			x y z\\
			0 0 0\\%0
			3 0 0\\%1
			-3 3 0\\%2
			0 -3 3\\%3
			 0 3 0\\%4
			-3 0 3\\%5
			3 -3 3\\%6
		};
		\draw[dashed, black!40] (axis cs: 0,0,0) -- (axis cs: 0,3,0);
		\draw[dashed, black!40] (axis cs: 0,0,0) -- (axis cs: 0,0,1);
		%\draw[dashed, gray] (axis cs: 0,0,0) -- (axis cs: 1,0,0);
		\draw[-, black] (axis cs: 0,0,0) -- (axis cs: 3,0,0);
		%\draw[dashed, gray] (axis cs: 0,0,0) -- (axis cs: 0,1,0);
		\draw[-, black] (axis cs: 0,0,1) -- (axis cs: 0,0,3);
	\end{axis}
\end{tikzpicture}

%% file: Bilder/FirstElemCone3D.tex
\begin{tikzpicture}[scale=0.65]
	\begin{axis}[
		view/h=30,
		view/v=30,
		xmin=-3, xmax=3.2,
		ymin=-3, ymax=3.2,
		zmin=-3, zmax=3.2,
		colormap={blau}{color=(blue!20) color=(blue!40) color=(blue!60)}
		]
		\addplot3[
		opacity=0.6,
		table/row sep=\\,
		patch,
		patch type=polygon,
		vertex count=4,
		patch table with point meta={%
			% pt1 pt2 pt3 pt4 pt5 cdata
			0 1 4 2 0\\
			0 3 5 2 1\\
			0 1 6 3 2\\
		}]
		table {
			x y z\\
			0 0 0\\%0
			-3 3 0\\%1
			0 -3 3\\%2
			0 0 -3\\%3
			-3 0 3\\%4
			0 -3 0\\%5
			-3 3 -3\\%6
		};
		\draw[dashed, black!20] (axis cs: 0,0,-2) -- (axis cs: 0,0,0);
		\draw[dashed, black] (axis cs: 0,-3,0) -- (axis cs: 0,0,0);
		%\draw[dashed, gray] (axis cs: 0,0,0) -- (axis cs: 1,0,0);
		\draw[-, black] (axis cs: -3,0,0) -- (axis cs: 0,0,0);
		%\draw[dashed, gray] (axis cs: 0,0,0) -- (axis cs: 0,1,0);
		%\draw[-, black] (axis cs: 0,1,0) -- (axis cs: 0,3.5,0);
	\end{axis}
\end{tikzpicture}

%% file: Bilder/LastElemDualCone3D.tex
\begin{tikzpicture}[scale=0.65]
	\begin{axis}[
		view/h=30,
		view/v=30,
		xmin=-3, xmax=3.2,
		ymin=-3, ymax=3.2,
		zmin=-3, zmax=3.2,
		colormap={blau}{color=(blue!20) color=(blue!40) color=(blue!60)}
		]
		\addplot3[
		opacity=0.6,
		table/row sep=\\,
		patch,
		patch type=polygon,
		vertex count=4,
		patch table with point meta={%
			% pt1 pt2 pt3 pt4 pt5 cdata
			0 1 4 2 2\\
			0 3 5 2 1\\
			0 1 6 3 0\\
		}]
		table {
			x y z\\
			0 0 0\\%0
			2 2 2\\%1
			0 2 2\\%2
			0 0 3\\%3
			2 4 4\\%4
			0 2 5\\%5
			2 2 5\\%6
		};
		%\draw[dashed, black!20] (axis cs: 0,0,-2) -- (axis cs: 0,0,0);
		%\draw[dashed, black] (axis cs: 0,-3,0) -- (axis cs: 0,0,0);
		%\draw[dashed, gray] (axis cs: 0,0,0) -- (axis cs: 1,0,0);
		%\draw[-, black] (axis cs: -3,0,0) -- (axis cs: 0,0,0);
		%\draw[dashed, gray] (axis cs: 0,0,0) -- (axis cs: 0,1,0);
		%\draw[-, black] (axis cs: 0,1,0) -- (axis cs: 0,3.5,0);
	\end{axis}
\end{tikzpicture}

%% file: Bilder/FirstElemDualCone3D.tex
\begin{tikzpicture}[scale=0.6]
	\begin{axis}[
		view/h=30,
		view/v=30,
		xmin=-3, xmax=3.2,
		ymin=-3, ymax=3.2,
		zmin=-3, zmax=3.2,
		colormap={blau}{color=(blue!20) color=(blue!40) color=(blue!60)}
		]
		\addplot3[
		opacity=0.6,
		table/row sep=\\,
		patch,
		patch type=polygon,
		vertex count=4,
		patch table with point meta={%
			% pt1 pt2 pt3 pt4 pt5 cdata
			0 1 4 2 2\\
			0 3 5 2 1\\
			0 1 6 3 0\\
		}]
		table {
			x y z\\
			0 0 0\\%0
			-3 0 0\\%1
			-2 -2 0\\%2
			-2 -2 -2\\%3
			-5 -2 0\\%4
			-4 -4 -2\\%5
			-5 -2 -2\\%6
		};
		%\draw[dashed, black!20] (axis cs: 0,0,-2) -- (axis cs: 0,0,0);
		%\draw[dashed, black] (axis cs: 0,-3,0) -- (axis cs: 0,0,0);
		%\draw[dashed, gray] (axis cs: 0,0,0) -- (axis cs: 1,0,0);
		\draw[-, black] (axis cs: 0,-3,0) -- (axis cs: 0,0,0);
		%\draw[dashed, gray] (axis cs: 0,0,0) -- (axis cs: 0,1,0);
		%\draw[-, black] (axis cs: 0,1,0) -- (axis cs: 0,3.5,0);
	\end{axis}
\end{tikzpicture}

%% file: Bilder/3K-COP-shadow.tex
%\documentclass {article}
%\usepackage{tikz}
%\usetikzlibrary{patterns}
%\usetikzlibrary{graphs,quotes}\usepackage{pgfplots}
%\pgfplotsset{compat=1.16,
%	every axis/.append style={
%		axis lines=center,
%		xlabel style={anchor=south west},
%		ylabel style={anchor=south west},
%		zlabel style={anchor=south west},
%		tick align=outside,}
%}
%\usepgfplotslibrary{patchplots}

%\begin{document}

\begin{tikzpicture}[scale=0.7]
	\clip (1,1.5) rectangle (6.6,4.5);
	\begin{axis}[
		view/h=-30,
		view/v=-30,
		 axis line style={draw=none},
		 xtick=\empty,
		 ytick=\empty,
		 ztick=\empty,
		xmin=-10, xmax=4.2,
		ymin=-10, ymax=4.2,
		zmin=-10, zmax=4.2,
		colormap={blau}{color=(blue!20) color=(blue!40) color=(blue!60)}
		]
		
		\addplot3[shading = darkGray2, draw=none]%axis, top color=black!20, bottom color=black!50, draw=none]
		coordinates{
			(4,4,4)
			(-10,4,4)
			(-10,-10,4)
			(4,-10,4)
		};
	
		\addplot3[shading = midGray, shading angle=-135, draw=none] %axis, right color=black!35, left color=black!10, shading angle=-135, draw=none]
		coordinates{
			(4,4,4)
			(4,-10,4)
			(4,-10,-10)
			(4,4,-10)
		};
	
		\addplot3[shading = lightGray2, shading angle=155, draw=none] %axis, left color=white, right color=black!25,shading angle=135, draw=none]
		coordinates{
			(4,4,4)
			(-10,4,4)
			(-10,4,-10)
			(4,4,-10)
		};

		\addplot3[
		table/row sep=\\,
		patch,
		patch type=polygon,
		vertex count=8,
		patch table with point meta={%
			% pt1 pt2 pt3 pt4 pt5 cdata
			0 1 2 3 0 0 0 0 0\\
			0 1 4 5 0 0 0 0 2\\
			0 6 7 8 9 10 11 5 1\\
			12 10 9 13 12 12 12 12 0\\
			12 10 11 14 12 12 12 12 2\\
			12 13 15 14 12 12 12 12 1\\
			16 8 7 17 16 16 16 16 0\\
			16 8 9 13 15 18 16 16 2\\
			16 17 19 18 16 16 16 16 1\\
		}]
		table {
			x y z\\
			3 1 0\\%0
			4 1 0\\%1
			4 -3 4\\%2
			3 -3 4\\%3
			4 4 0\\%4
			0 4 0\\%5
			3 -3 4\\%6
			2 -2 4\\ %7
			2 0 2\\%8
			1 1 2\\%9
			1 2 1\\%10
			-1 4 1\\%11
			0 2 1\\%12
			0 1 2\\%13
			-2 4 1\\%14
			-3 4 2\\%15
			0 0 2\\%16
			0 -2 4\\%17
			-4 4 2\\%18
			-6 4 4\\%19
		};
	
		\addplot3 [only marks, mark=*, mark size =1.5pt] coordinates{(3,1,0) };
		\node at (axis cs: 2.5 ,0.5 ,0.5) {$c^1$ };
		\addplot3 [only marks, mark=*, mark size =1.5pt] coordinates{(0,2,1) };
		\node at (axis cs: -0.5 ,1.5 ,1.3) {$c^2$ };
		\addplot3 [only marks, mark=*, mark size =1.5pt] coordinates{(0,0,2) };
		\node at (axis cs: -0.5 ,-0.5 ,2.5) {$c^3$ };
		\addplot3 [only marks, mark=*, mark size =1.5pt] coordinates{(1,0,2) };
		\node at (axis cs: 1 ,-0.5 ,2.8) {$c^4$ };
		
		%\draw[-, black] (axis cs: 0,-3,0) -- (axis cs: 0,5,0);
		%\draw[-, black] (axis cs: -3,0,0) -- (axis cs: 5,0,0);
		%\draw[-, black] (axis cs: 0,0,-3) -- (axis cs: 0,0,5);
		
		%\draw[-, black] (axis cs: -7,-7,-7) -- (axis cs: -7,4,-7);
		%\draw[-, black] (axis cs: -7,-7,-7) -- (axis cs: 4,-7,-7);
		%\draw[-, black] (axis cs: -7,-7,-7) -- (axis cs: -7,-7,4);
		%\draw[-, black] (axis cs: 4,-7,-7) -- (axis cs: 4,-7,4);
		%\draw[-, black] (axis cs: 4,-7,-7) -- (axis cs: 4,4,-7);
		%\draw[-, black] (axis cs: -7,4,-7) -- (axis cs: 4,4,-7);
		%\draw[-, black] (axis cs: -7,4,-7) -- (axis cs: -7,4,4);
		%\draw[-, black] (axis cs: -7,4,4) -- (axis cs: -7,-7,4);
		%\draw[dashed, gray] (axis cs: -6,4,4) -- (axis cs: -7,4,4);
		%\draw[dashed, gray] (axis cs: 4,-3,4) -- (axis cs: 4,-7,4);
		%\draw[dashed, gray] (axis cs: 4,4,0) -- (axis cs: 4,4,-7);
		%\draw[-, black] (axis cs: -7,-7,4) -- (axis cs: 4,-7,4);
	\end{axis}
\end{tikzpicture}

%\end{document}

%% file: Bilder/3K-POP-shadow.tex
%\documentclass {article}
%\usepackage{tikz}
%\usetikzlibrary{patterns}
%\usetikzlibrary{graphs,quotes}\usepackage{pgfplots}
%\pgfplotsset{compat=1.16,
%	every axis/.append style={
%		axis lines=center,
%		xlabel style={anchor=south west},
%		ylabel style={anchor=south west},
%		zlabel style={anchor=south west},
%		tick align=outside,}
%}
%\usepgfplotslibrary{patchplots}

%\begin{document}

\begin{tikzpicture}[scale=0.7]
	\clip (1,1.5) rectangle (6.6,4.5);
	\begin{axis}[
		view/h=-30,
		view/v=-30,
		 axis line style={draw=none},
		 xtick=\empty,
		 ytick=\empty,
		 ztick=\empty,
		xmin=-6, xmax=5.2,
		ymin=-6, ymax=5.2,
		zmin=-6, zmax=5.2,
		colormap={blau}{color=(blue!20) color=(blue!40) color=(blue!60)}
		]
		
		\addplot3[shading = darkGray, shading angle=5,draw=none]
		coordinates{
			(5,5,5)
			(-10,5,5)
			(-10,-10,5)
			(5,-10,5)
		};
		
		\addplot3[shading = midGray, shading angle=-110, draw=none]
		coordinates{
			(5,5,5)
			(5,-10,5)
			(5,-10,-10)
			(5,5,-10)
		};
		
		\addplot3[shading = lightGray,shading angle=110, draw=none]
		coordinates{
			(5,5,5)
			(-10,5,5)
			(-10,5,-10)
			(5,5,-10)
		};
		
		\addplot3[
		table/row sep=\\,
		patch,
		patch type=polygon,
		vertex count=8,
		patch table with point meta={%
			% pt1 pt2 pt3 pt4 pt5 cdata
			0 1 2 3 0 0 0 0 0\\
			0 1 4 5 0 0 0 0 2\\
			0 6 7 8 9 10 11 5 1\\
			12 10 9 13 12 12 12 12 0\\
			12 10 11 14 12 12 12 12 2\\
			12 13 15 14 12 12 12 12 1\\
			16 8 7 17 16 16 16 16 0\\
			16 8 9 13 15 18 16 16 2\\
			16 17 19 18 16 16 16 16 1\\
		}]
		table {
			x y z\\
			4 1 0\\%0
			5 1 0\\%1
			5 1 5\\%2
			4 1 5\\%3
			5 5 0\\%4
			4 5 0\\%5
			4 1 5\\%6
			4 2 5\\ %7
			4 2 2\\%8
			4 3 2\\%9
			4 3 1\\%10
			4 5 1\\%11
			3 3 1\\%12
			3 3 2\\%13
			3 5 1\\%14
			3 5 2\\%15
			2 2 2\\%16
			2 2 5\\%17
			2 5 2\\%18
			2 5 5\\%19
		};
	
		\addplot3 [only marks, mark=*, mark size =1.5pt] coordinates{(4,1,0) };
		\node at (axis cs: 3.2 ,0.6 ,0.2) {$\tilde{c}^1$ };
		\addplot3 [only marks, mark=*, mark size =1.5pt] coordinates{(3, 3, 1) };
		\node at (axis cs: 2.3 ,2.6 ,1.2) {$\tilde{c}^2$ };
		\addplot3 [only marks, mark=*, mark size =1.5pt] coordinates{(2,2,2) };
		\node at (axis cs: 1.2 ,1.6 ,2.1) {$\tilde{c}^3$ };
		\addplot3 [only marks, mark=*, mark size =1.5pt] coordinates{(3,2,2) };
		\node at (axis cs: 3 ,1.6 ,2.7) {$\tilde{c}^4$ };
		
		%\draw[-, black] (axis cs: 0,-3,0) -- (axis cs: 0,5,0);
		%\draw[-, black] (axis cs: -3,0,0) -- (axis cs: 5,0,0);
		%\draw[-, black] (axis cs: 0,0,-3) -- (axis cs: 0,0,5);
		
		%\draw[-, black] (axis cs: -3,-3,-3) -- (axis cs: -3,5,-3);
		%\draw[-, black] (axis cs: -3,-3,-3) -- (axis cs: 5,-3,-3);
		%\draw[-, black] (axis cs: -3,-3,-3) -- (axis cs: -3,-3,5);
		%\draw[-, black] (axis cs: 5,-3,-3) -- (axis cs: 5,-3,5);
		%\draw[-, black] (axis cs: 5,-3,-3) -- (axis cs: 5,5,-3);
		%\draw[-, black] (axis cs: -3,5,-3) -- (axis cs: 5,5,-3);
		%\draw[-, black] (axis cs: -3,5,-3) -- (axis cs: -3,5,5);
		%\draw[-, black] (axis cs: -3,5,5) -- (axis cs: -3,-3,5);
		%\draw[dashed, gray] (axis cs: 2,5,5) -- (axis cs: -3,5,5);
		%\draw[dashed, gray] (axis cs: 5,1,5) -- (axis cs: 5,-3,5);
		%\draw[dashed, gray] (axis cs: 5,5,0) -- (axis cs: 5,5,-3);
		%\draw[-, black] (axis cs: -3,-3,5) -- (axis cs: 5,-3,5);
	\end{axis}

\end{tikzpicture}

%\end{document}

%% file: Bilder/WeightSpaceDec.tex
\begin{tikzpicture}[scale=4]
	% WSD
	\fill[fill=black!20,fill opacity=1] (0.33333333333333,0.33333333333333) -- (0,1) -- (0,0)  -- cycle;
	\fill[fill=black!40,fill opacity=1] (0.33333333333333,0.33333333333333) -- (0,1) -- (1,0)  -- cycle;
	\fill[fill=black!60,fill opacity=1] (0.33333333333333,0.33333333333333) -- (1,0) -- (0,0)  -- cycle;
	\draw (0,0) -- (0.33333333333333,0.33333333333333);
	\draw (1,0) -- (0.33333333333333,0.33333333333333);
	\draw (0,1) -- (0.33333333333333,0.33333333333333);
	\draw[dashed] (1,0) -- (0,1);
	
	% axes
	\draw[dashed,thick] (0,0) -- (1,0);
	\draw[->,thick] (1,0) -- (1.1,0) node[right] {$\lambda_1$};
	\draw[dashed,thick] (0,0) -- (0,1);
	\draw[->,thick] (0,1) -- (0,1.1) node[above] {$\lambda_2$};
	\draw (-0.02,0.16666666666666) -- (0.02,0.16666666666666) node[left=5pt] {$\frac{1}{6}$};
	\draw (-0.02,0.33333333333333) -- (0.02,0.33333333333333) node[left=5pt] {$\frac{1}{3}$};
	\draw (-0.02,0.5) -- (0.02,0.5) node[left=5pt] {$\frac{1}{2}$};
	\draw (-0.02,0.66666666666666) -- (0.02,0.66666666666666) node[left=5pt] {$\frac{2}{3}$};
	\draw (-0.02,0.83333333333333) -- (0.02,0.83333333333333) node[left=5pt] {$\frac{5}{6}$};
	\draw (-0.02,1) -- (0.02,1) node[left=5pt] {$1$};
	\draw (0.16666666666666,-0.02) -- (0.16666666666666,0.02) node[below=5pt] {$\frac{1}{6}$};
	\draw (0.33333333333333,-0.02) -- (0.33333333333333,0.02) node[below=5pt] {$\frac{1}{3}$};
	\draw (0.5,-0.02) -- (0.5,0.02) node[below=5pt] {$\frac{1}{2}$};
	\draw (0.66666666666666,-0.02) -- (0.66666666666666,0.02) node[below=5pt] {$\frac{2}{3}$};
	\draw (0.83333333333333,-0.02) -- (0.83333333333333,0.02) node[below=5pt] {$\frac{5}{6}$};
	\draw (1,-0.02) -- (1,0.02) node[below=5pt] {$1$};
\end{tikzpicture}

%% file: Bilder/OrdinalSpaceDec.tex
\begin{tikzpicture}[scale=4]
	% WSD
	\fill[fill=black!20,fill opacity=1] (0.16666666666666,0.33333333333333) -- (0,0.5) -- (0,0)  -- cycle;
	\fill[fill=black!40,fill opacity=1] (0.16666666666666,0.33333333333333) -- (0,0.5) -- (0.33333333333333,0.33333333333333)  -- cycle;
	\fill[fill=black!60,fill opacity=1] (0.16666666666666,0.33333333333333) -- (0.33333333333333,0.33333333333333) -- (0,0)  -- cycle;
	\draw (0,0) -- (0.16666666666666,0.33333333333333);
	\draw (0.33333333333333,0.33333333333333) -- (0.16666666666666,0.33333333333333);
	\draw (0,0.5) -- (0.16666666666666,0.33333333333333);
	\draw[dashed] (0.33333333333333,0.33333333333333) -- (0,0.5);
	\draw[dashed] (0.33333333333333,0.33333333333333) -- (0,0);
	
	% axes
	\draw[->,thick] (0,0) -- (1.1,0) node[right] {$\mu_1$};
	\draw[dashed,thick] (0,0) -- (0,0.5);
	\draw[->,thick] (0,0.5) -- (0,1.1) node[above] {$\mu_2$};
	\draw (-0.02,0.16666666666666) -- (0.02,0.16666666666666) node[left=5pt] {$\frac{1}{6}$};
	\draw (-0.02,0.33333333333333) -- (0.02,0.33333333333333) node[left=5pt] {$\frac{1}{3}$};
	\draw (-0.02,0.5) -- (0.02,0.5) node[left=5pt] {$\frac{1}{2}$};
	\draw (-0.02,0.66666666666666) -- (0.02,0.66666666666666) node[left=5pt] {$\frac{2}{3}$};
	\draw (-0.02,0.83333333333333) -- (0.02,0.83333333333333) node[left=5pt] {$\frac{5}{6}$};
	\draw (-0.02,1) -- (0.02,1) node[left=5pt] {$1$};
	\draw (0.16666666666666,-0.02) -- (0.16666666666666,0.02) node[below=5pt] {$\frac{1}{6}$};
	\draw (0.33333333333333,-0.02) -- (0.33333333333333,0.02) node[below=5pt] {$\frac{1}{3}$};
	\draw (0.5,-0.02) -- (0.5,0.02) node[below=5pt] {$\frac{1}{2}$};
	\draw (0.66666666666666,-0.02) -- (0.66666666666666,0.02) node[below=5pt] {$\frac{2}{3}$};
	\draw (0.83333333333333,-0.02) -- (0.83333333333333,0.02) node[below=5pt] {$\frac{5}{6}$};
	\draw (1,-0.02) -- (1,0.02) node[below=5pt] {$1$};
\end{tikzpicture}

%% file: Bilder/2K+w-COP.tex
\begin{tikzpicture}[scale=0.95,every node/.style={font=\footnotesize}]
	\clip (1,1.5) rectangle (6.6,4.5);
	\begin{axis}[
		view/h=-30,
		view/v=-30,
		 axis line style={draw=none},
		 xtick=\empty,
		 ytick=\empty,
		 ztick=\empty,
		xmin=-8, xmax=5.2,
		ymin=-8, ymax=5.2,
		zmin=-8, zmax=5.2,
		colormap={blau}{color=(blue!20) color=(blue!40) color=(blue!60)}
		]
		
		\addplot3[shading = darkGray, shading angle=25,draw=none]
		coordinates{
			(5,5,5)
			(-10,5,5)
			(-10,-10,5)
			(5,-10,5)
		};
		
		\addplot3[shading = midGray, shading angle=-120, draw=none]
		coordinates{
			(5,5,5)
			(5,-10,5)
			(5,-10,-10)
			(5,5,-10)
		};
		
		\addplot3[shading = lightGray,shading angle=110, draw=none]
		coordinates{
			(5,5,5)
			(-10,5,5)
			(-10,5,-10)
			(5,5,-10)
		};
		
		\addplot3[
		table/row sep=\\,
		patch,
		patch type=polygon,
		vertex count=8,
		patch table with point meta={%
			% pt1 pt2 pt3 pt4 pt5 cdata
			0 1 2 3 0 0 0 0 0\\
			0 1 4 5 0 0 0 0 2\\
			0 6 7 8 9 10 11 5 1\\
			12 10 9 13 12 12 12 12 0\\
			12 10 11 14 12 12 12 12 2\\
			12 13 15 14 12 12 12 12 1\\
			16 8 7 17 16 16 16 16 0\\
			16 8 9 13 15 18 16 16 2\\
			16 17 19 18 16 16 16 16 1\\
		}]
		table {
			x y z\\
			4 1 0\\%0
			5 1 0\\%1
			5 -4 5\\%2
			4 -4 5\\%3
			5 5 0\\%4
			4 5 0\\%5
			4 -4 5\\%6
			4 -3 5\\ %7
			4 0 2\\%8
			4 1 2\\%9
			4 2 1\\%10
			4 5 1\\%11
			3 2 1\\%12
			3 1 2\\%13
			3 5 1\\%14
			3 5 2\\%15
			2 0 2\\%16
			2 -3 5\\%17
			2 5 2\\%18
			2 5 5\\%19
		};
	
		\addplot3 [only marks, mark=*] coordinates{(4,1,0) };
		\node at (axis cs: 3.5 ,0.5 ,0.5) {$v^1$ };
		\addplot3 [only marks, mark=*] coordinates{(3,2,1) };
		\node at (axis cs: 2.5 ,1.5 ,1.3) {$v^2$ };
		\addplot3 [only marks, mark=*] coordinates{(2,0,2) };
		\node at (axis cs: 1.5 ,-0.5 ,2.5) {$v^3$ };
		\addplot3 [only marks, mark=*] coordinates{(3,0,2) };
		\node at (axis cs: 3 ,-0.5 ,2.7) {$v^4$ };
		
		%\draw[-, black] (axis cs: 0,-3,0) -- (axis cs: 0,5,0);
		%\draw[-, black] (axis cs: -3,0,0) -- (axis cs: 5,0,0);
		%\draw[-, black] (axis cs: 0,0,-3) -- (axis cs: 0,0,5);
		
		%\draw[-, black] (axis cs: -5,-5,-5) -- (axis cs: -5,5,-5);
		%\draw[-, black] (axis cs: -5,-5,-5) -- (axis cs: 5,-5,-5);
		%\draw[-, black] (axis cs: -5,-5,-5) -- (axis cs: -5,-5,5);
		%\draw[-, black] (axis cs: 5,-5,-5) -- (axis cs: 5,-5,5);
		%\draw[-, black] (axis cs: 5,-5,-5) -- (axis cs: 5,5,-5);
		%\draw[-, black] (axis cs: -5,5,-5) -- (axis cs: 5,5,-5);
		%\draw[-, black] (axis cs: -5,5,-5) -- (axis cs: -5,5,5);
		%\draw[-, black] (axis cs: -5,5,5) -- (axis cs: -5,-5,5);
		%\draw[dashed, gray] (axis cs: 2,5,5) -- (axis cs: -5,5,5);
		%\draw[dashed, gray] (axis cs: 5,-4,5) -- (axis cs: 5,-5,5);
		%\draw[dashed, gray] (axis cs: 5,5,0) -- (axis cs: 5,5,-5);
		%\draw[-, black] (axis cs: -5,-5,5) -- (axis cs: 5,-5,5);
	\end{axis}
\end{tikzpicture}

%% file: Bilder/2K+w-POP.tex
\begin{tikzpicture}[scale=0.95,every node/.style={font=\footnotesize}]
	\clip (1,1.5) rectangle (6.6,4.5);
	\begin{axis}[
		view/h=-30,
		view/v=-30,
		 axis line style={draw=none},
		 xtick=\empty,
		 ytick=\empty,
		 ztick=\empty,
		xmin=-8, xmax=5.2,
		ymin=-8, ymax=5.2,
		zmin=-8, zmax=5.2,
		colormap={blau}{color=(blue!20) color=(blue!40) color=(blue!60)}
		]

		\addplot3[shading = darkGray, shading angle=5,draw=none]
		coordinates{
			(5,5,5)
			(-10,5,5)
			(-10,-10,5)
			(5,-10,5)
		};
		
		\addplot3[shading = midGray, shading angle=-110, draw=none]
		coordinates{
			(5,5,5)
			(5,-10,5)
			(5,-10,-10)
			(5,5,-10)
		};
		
		\addplot3[shading = lightGray,shading angle=110, draw=none]
		coordinates{
			(5,5,5)
			(-10,5,5)
			(-10,5,-10)
			(5,5,-10)
		};
		
		\addplot3[
		table/row sep=\\,
		patch,
		patch type=polygon,
		vertex count=8,
		patch table with point meta={%
			% pt1 pt2 pt3 pt4 pt5 cdata
			0 1 2 3 0 0 0 0 0\\
			0 1 4 5 0 0 0 0 2\\
			0 6 7 8 9 10 11 5 1\\
			12 10 9 13 12 12 12 12 0\\
			12 10 11 14 12 12 12 12 2\\
			12 13 15 14 12 12 12 12 1\\
			16 8 7 17 16 16 16 16 0\\
			16 8 9 13 15 18 16 16 2\\
			16 17 19 18 16 16 16 16 1\\
		}]
		table {
			x y z\\
			4 1 0\\%0
			5 1 0\\%1
			5 1 5\\%2
			4 1 5\\%3
			5 5 0\\%4
			4 5 0\\%5
			4 1 5\\%6
			4 2 5\\ %7
			4 2 2\\%8
			4 3 2\\%9
			4 3 1\\%10
			4 5 1\\%11
			3 3 1\\%12
			3 3 2\\%13
			3 5 1\\%14
			3 5 2\\%15
			2 2 2\\%16
			2 2 5\\%17
			2 5 2\\%18
			2 5 5\\%19
		};
	
		\addplot3 [only marks, mark=*] coordinates{(4,1,0) };
		\node at (axis cs: 3.5 ,0.6 ,0.2) {$\tilde{v}^1$ };
		\addplot3 [only marks, mark=*] coordinates{(3, 3, 1) };
		\node at (axis cs: 2.5 ,2.6 ,1.2) {$\tilde{v}^2$ };
		\addplot3 [only marks, mark=*] coordinates{(2,2,2) };
		\node at (axis cs: 1.3 ,1.6 ,2.1) {$\tilde{v}^3$ };
		\addplot3 [only marks, mark=*] coordinates{(3,2,2) };
		\node at (axis cs: 3.3 ,1.6 ,2.9) {$\tilde{v}^4$ };
		
		%\draw[-, black] (axis cs: 0,-3,0) -- (axis cs: 0,5,0);
		%\draw[-, black] (axis cs: -3,0,0) -- (axis cs: 5,0,0);
		%\draw[-, black] (axis cs: 0,0,-3) -- (axis cs: 0,0,5);
		
		%\draw[-, black] (axis cs: -3,-3,-3) -- (axis cs: -3,5,-3);
		%\draw[-, black] (axis cs: -3,-3,-3) -- (axis cs: 5,-3,-3);
		%\draw[-, black] (axis cs: -3,-3,-3) -- (axis cs: -3,-3,5);
		%\draw[-, black] (axis cs: 5,-3,-3) -- (axis cs: 5,-3,5);
		%\draw[-, black] (axis cs: 5,-3,-3) -- (axis cs: 5,5,-3);
		%\draw[-, black] (axis cs: -3,5,-3) -- (axis cs: 5,5,-3);
		%\draw[-, black] (axis cs: -3,5,-3) -- (axis cs: -3,5,5);
		%\draw[-, black] (axis cs: -3,5,5) -- (axis cs: -3,-3,5);
		%\draw[dashed, gray] (axis cs: 2,5,5) -- (axis cs: -3,5,5);
		%\draw[dashed, gray] (axis cs: 5,1,5) -- (axis cs: 5,-3,5);
		%\draw[dashed, gray] (axis cs: 5,5,0) -- (axis cs: 5,5,-3);
		%\draw[-, black] (axis cs: -3,-3,5) -- (axis cs: 5,-3,5);
	\end{axis}
\end{tikzpicture}

%% file: Bilder/graphExtension.tex
\begin{tikzpicture}[xscale=1.5, yscale=1.7,every edge quotes/.append style={font=\footnotesize}]%[shorten >=1pt,->,draw=black!75, node distance=\layersep,scale=0.28]
\node[draw,circle] (1) at (0,0)[] {$s$};
\node[draw,circle] (2) at (1,0.5) {};
\node[draw,circle] (3) at (2,0.5) {};
\node[draw,circle] (4) at (3,0) {$t$};
\node[draw,circle] (5) at (1,-0.5) {};
\node[draw,circle] (6) at (2,-0.5) {};

%\node (BotLeft) at (-2,-6) {};
%\node (BotLeft2) at (-1.75,-4.75) {};

%%graph G
\graph {
	(1) ->["$w(e_1)=1$",-latex,thick,draw=black!40!red] (2);
	(2) ->["$w(e_2)=1$",-latex,thick,draw=black!40!red] (3);
	(3) ->["$w(e_3)=8$",-latex,densely dotted,thick,draw=black!50!green] (4);
	(1) ->["$w(e_4)=6$",-latex,thick,swap,draw=black!40!red] (5);
	(5) ->["$w(e_5)=2$",-latex,densely dotted,thick,swap,draw=black!50!green] (6);
	(6) ->["$w(e_6)=2$",-latex,densely dotted,swap,thick,draw=black!50!green] (4);
};
\end{tikzpicture}